\documentclass[smallextended]{svjour3}
\usepackage{pifont}
\usepackage{epsfig,color,graphicx,subfigure,amsmath,amssymb,multirow,lscape}
\usepackage[colorlinks,linkcolor=blue,citecolor=blue]{hyperref}
\usepackage{cite}
\smartqed
\usepackage{graphicx}
\usepackage{amsmath,amssymb}
\journalname{COAP}

\newtheorem{thm}{Theorem}
\newtheorem{lem}{Lemma}
\newtheorem{rem}{Remark}
\newtheorem{fact}{Fact}
\newtheorem{algo}{Algorithm}
\newtheorem{prob}{Problem}

\def\be{\begin{eqnarray}}
\def\ee{\end{eqnarray}}
\def\ben{\begin{eqnarray*}}
\def\een{\end{eqnarray*}}
\def\ba{\begin{array}}
\def\ea{\end{array}}
\def\proof {\pn {Proof.} }
\def\endproof{\hfill $\Box$ \vskip .5cm}

\def\cH{{\mathcal H}}
\def\cO{{\mathcal O}}
\def\cS{{\mathcal S}}

\def\bR{{\mathbb R}}
\def\bN{{\mathbb N}}
\def\pn{\par\smallskip\noindent}

\def\prox{{\rm Prox}}

\def\proof {\pn {Proof.} }

\newcommand{\blue}[1]{\begin{color}{blue}#1\end{color}}

\DeclareMathOperator*{\dom}{dom}
\DeclareMathOperator*{\argmin}{argmin}

\begin{document}

\title{Proximal extrapolated gradient methods with prediction and correction for monotone variational inequalities }
\titlerunning{Proximal extrapolated gradient method}

\author{Xiaokai Chang$^{1,2}$   \and  Sanyang Liu$^1$ \and Jianchao Bai$^3$ \and Jun Yang$^4$ }

\institute{\ding {41} Xiaokai Chang \at xkchang@lut.cn \\
\ding {41} Jianchao Bai \at bjc1987@163.com
           \and
1 \ \   School of Science, Lanzhou University of Technology, Lanzhou, P. R. China. \\
2  \ \  School of Mathematics and Statistics, Xidian University, Xi'an, P. R. China. \\
3  \ \  Department of Applied Mathematics, Northwestern Polytechnical University, Xi'an, P. R. China.\\
4  \ \  School of Mathematics and Information Science, Xianyang Normal University, Xianyang, P. R. China.}

\date{Received: date / Accepted: date}

\maketitle

\begin{abstract}
An efficient proximal-gradient-based method, called proximal extrapolated gradient method, is designed for solving monotone variational inequality in Hilbert space.
The proposed method extends the acceptable range of parameters to obtain larger step sizes. The step size is predicted based a local information of the operator and corrected by linesearch procedures to satisfy a very weak condition, which is even weaker than the boundedness of sequence generated and always holds when the operator is the gradient of a convex function.  We establish  its convergence and    ergodic convergence rate in theory under the larger range of parameters.  Furthermore, we improve numerical efficiency by employing the proposed method with non-monotonic step size, and obtain the upper bound of the parameter relating to step size by an extremely simple example.
Related  numerical experiments illustrate the improvements in efficiency from the larger step size.
\end{abstract}

\section{Introduction}
\label{sec_introduction}

Let $\mathcal{H}$ be a real Hilbert space equipped with inner product $\langle \cdot , \cdot\rangle   $ and its induced norm $\parallel\cdot\parallel $. We consider the variational inequality problem:
\be\label{primal_problem0}
\mbox{find} ~~ x^*\in \cH   ~~ \textrm{s.t.} ~\langle F(x^*),y-x^*\rangle+g(y)-g(x^*)\geq 0,~ \forall y\in \cH,
\ee
where $F :\cH\rightarrow\cH$ is an operator and $g :\cH\rightarrow]-\infty, +\infty]$ is a proper lower semicontinuous convex function. We use $\dom g$ to represent the domain of $g$, defined by $\dom g := \{x\in\cH: g(x)<+\infty\}$. For a continuously differentiable and convex function $f :\cH\rightarrow]-\infty, +\infty[$ with its gradient denoted by $\nabla f=F$, then problem (\ref{primal_problem0}) is equivalent to
\be\label{two-block}
\min_{x\in\cH} f(x)+g(x).
\ee
Let $C$ be a closed and convex subset of $\cH$. Let $l_{C}$ be the indicator function of the set $C$, that is, $l_{C}(x)=0$ if $x\in C$ and $\infty$ otherwise. When $g(x) = l_{C}(x)$, variational inequality (\ref{primal_problem0}) reduces to
\be\label{primal_problem}
\mbox{find} ~~ x^*\in C   ~~ \textrm{s.t.} ~\langle F(x^*),y-x^*\rangle\geq 0,~ \forall y\in \cH.
\ee
Problem (\ref{primal_problem0}) and its special cases (\ref{two-block}) and (\ref{primal_problem}) have wide applications    in   disciplines  including mechanics, signal and image processing, and economics \cite{app1,app2,app3,app4,statistical_learning,11.}, to cite a few. Throughout the paper,  the solution set $\cS$ of problem (\ref{primal_problem0}) is assumed to be nonempty, and the following assumptions hold: \par
$\textbf{(A1)}$ \   $F$ is monotone, i.e.,
  \begin{center}
 $\langle F(x)-F(y),x-y\rangle\geq0,
  \ \ \forall x, y\in \mathcal{H};$
  \end{center}\par
$\textbf{(A2)}$ \  $F$ is $L$-Lipschitz continuous ($L>0$), that is,
\ben
\| F(x)-F(y) \|\leq L\| x-y \|,\ \ \forall x, y\in \cH;
\een\par
$\textbf{(A3)}$ \  $g|\dom g$ is a continuous function.

Many efficient methods have been proposed for solving the problem (\ref{primal_problem0}) and its special cases, for instance,  alternating direction method of multipliers (ADMM) \cite{statistical_learning,PC-ADMM,He_ADMM-based,ADMM}, extragradient method \cite{6.,extragradient,extragradient-type,13.}, proximal (projected) gradient method  \cite{8.,FBS,FBS_P,FBS2015,modified-FB,New_properties} and its accelerated version \cite{FISTA-CD,Nesterov1983}. Here, we would  concentrate on the most simple case of these approaches: forward-backward splitting (FBS) method. Under the assumption that $F$ is $L$-Lipschitz continuous, the iterative scheme of the classical FBS method for problem (\ref{primal_problem0}) reads
\be\label{classic_pg}
x_{n+1}=\mbox{prox}_{\lambda g}(x_n-\lambda F(x_n)),
\ee
where $\lambda$ is some positive number and can be viewed as a step size of the forward step, and the proximal operator $\mbox{prox}_{\lambda g} : \cH\rightarrow\cH$ is defined in Section \ref{sec_preliminarries}.

To establish   convergence of the iteration (\ref{classic_pg}), it often requires the restrictive assumptions that $F$ is $L$-Lipschitz continuous, strongly (or inverse strongly) monotone with $\lambda\in]0,\frac{2}{L}[$. To overcome this drawback, Korpelevich \cite{extragradient} and Antipin \cite{13.} proposed the following extragradient method for   (\ref{primal_problem}) with two-step projection procedures
\ben
y_n=P_C(x_n-\lambda_n F(x_n)),\ \ \ \ \  x_{n+1}=P_C(x_n-\lambda_n F(y_n)),
\een
where $P_C:\mathcal{H}\rightarrow C$ denotes the (metric) projection onto $C$, $\lambda_n$ is any positive sequence verifying $\lambda_n\in[l, u]$ for some values $l,u\in]0, \frac{1}{L}[$. The extragradient method has received great attentions and has been improved in various ways \cite{5.,M-extra,Non-Lip,Low-cost,9.}, including linesearch procedures or/and avoiding   Lipschitz-continuity assumption, decreasing a number of metric projections, etc. For instance, Censor, Gibali and Reich \cite{5.} introduced
\be\label{Censor}
\left.
\ba{l}
y_n=P_C(x_n-\lambda F(x_n)),\\
T_n=\{w\in\cH | \langle x_n-\lambda F(x_n)-y_n,~w-y_n\rangle\leq0\},\\
x_{n+1}=P_{T_n}(x_n-\lambda F(y_n)),
\ea\right\}
\ee
where the step size satisfies $\lambda \in ]0,\frac{1}{L}[$. Since the second projection $P_{T_n}$ in (\ref{Censor}) can be found in a closed form, this method is more applicable when a projection onto the closed convex set $C$ is a nontrivial problem. For  a more general problem (\ref{primal_problem0}), Tseng \cite{modified-FB} modified the iteration (\ref{classic_pg}) and proposed the following forward-backward-forward (FBF) method involving one proximal operator and two values of $F$ per iteration:
\ben
y_n=\mbox{prox}_{\lambda g}(x_n-\lambda F(x_n)), \ \ \ \ x_{n+1}=y_n+\lambda( F(x_n)-F(y_n)),
\een
where $\lambda \in ]0,\frac{1}{L}[$. Since then, Tseng's method has attracted a lot of interests  due to its simplicity and
generality, see \cite{FB-Tseng,Tseng,inertial-FBF} for more details.

In the literature,  the  inertial extrapolation has been conducted to accelerated proximal gradient methods in
the spirit of Nesterov's extrapolation techniques \cite{Nesterov1983,Nesterov2004}, whose basic idea is to make full use of historical information at each iteration. A typical scheme of the proximal gradient method with extrapolation for solving (\ref{primal_problem0}) is
\be\label{e_pg}
x_{n+1}=\prox_{\lambda g}(x_n-\lambda F(y_n)),~~y_{n+1}= x_{n+1}+\delta_n(x_{n+1}-x_n),
\ee
where $\delta_n>0$. Recently, using a fixed parameter $\delta=1$ in (\ref{e_pg}), Malitsky \cite{9.} introduced the iteration
\ben
x_{n+1}=P_C(x_n-\lambda F(2x_n-x_{n-1} )),~~\lambda \in ]0,(\sqrt2-1)/L[,
\een
for solving (\ref{primal_problem}). However, the step size ($\lambda$ or $\lambda_n$) requires the  information of the Lipschitz constant $L$, which is a main drawback of the algorithms introduced above.  In fact, these algorithms with a large value of $L$ can lead to very small step size, which may give rise to a slow convergent algorithm \cite{10.}. To obtain a proper step size, Armijo-type line search and outer approximation techniques were involved in \cite{Khobotov,search-strategy,Solodov,extragradient-type}. Due to the extra proximal operator as well as the  evaluations of $F$,  these algorithms will   be computationally expensive   when   proximal operator or $F$ is hard to compute and somewhat expensive.

For getting a proper step without using the Lipschitz constant $L$, Malitsky \cite{9.} introduced an efficient method whose main updates are
\begin{equation} \label{eq:5}
  \left\{ \begin{aligned}\begin{array}{llll}
          \mbox{Choose}\ \  x_0=y_0\in \mathcal{H},\  \lambda_0>0, \  \alpha \in ]0, {\sqrt2-1}[\\
          \mbox{Choose}\ \lambda_n \ \  \ s.t. \ \lambda_n \parallel Fy_n-Fy_{n-1} \parallel \leq \alpha \parallel y_n-y_{n-1} \parallel      \\
          x_{n+1}=P_C(x_n-\lambda_n F(y_n )) \\
          y_{n+1}= 2x_{n+1}-x_n .     \end{array}        \end{aligned} \right.
                            \end{equation}
By updating the step size $\lambda_n$ given by a specific procedure according the progress of algorithm, a weak convergence result was proved, but this process involves the computation of additional projections onto $C$. Later, Mainge and Gobinddass \cite{10.} introduced a more general framework:
\ben
\theta_n=\frac{\lambda_n}{\delta \lambda_{n-1}},~~~y_n=x_n+\theta_n(x_n-x_{n-1}),  ~~~ x_{n+1}=P_C(x_n-\lambda_n F(y_n )),
\een
where the step size $\lambda_n$ needs to satisfy
many inequality constraints and can be obtained by linesearch procedure, see
\cite[Section 3.1 and Section 3.2.2]{10.}. Based on the scheme (\ref{eq:5}), local information of the operator and some linesearch procedures, Malitsky \cite{Proximal-extrapolated}  proposed simpler schemes  which do not require Lipschitz continuity of the operator. Furthermore,  the involved linesearch procedure  doesn't need extra prox or projection and it can be applied to a more general problem (\ref{primal_problem0}). By overcoming the estimation of $L$ and linesearch procedure for the scheme (\ref{eq:5}),   Yang and Liu \cite{yang} proposed an extragradient method with lower computational complexity but   nonincreasing step sizes. The important parameter $\alpha$ relating to the step size $\lambda_n$ was restricted on $\alpha \in ]0, \frac{\sqrt{2}-1}{\delta}[$ with $\delta\in]1,+\infty[$ in \cite{yang} and $\alpha \in ]0, \sqrt{2}-1[$ with variable $\delta_n$ from linesearch in \cite{9.} for guaranteeing the convergence.

The aim of this paper is to propose a proximal gradient algorithm with larger step size, extend the range of $\delta$ to that is less than or equal to 1, and then improve the range of $\alpha$. Our proposed methods do not require Lipschitz constant, and its step size is predicted by using two previous iterates, and corrected by linesearch to satisfy a very weak condition, which always holds when $F=\nabla f$ for a convex function $f$.  Specifically, by the aid of the vital inequalities in convergence's proof we first   introduce a function $\kappa(\delta)$  defined as
\be \label{kappa}
\kappa(\delta):=\max\limits_{\varepsilon_1>0,\varepsilon_2>0}\min \left\{\frac{\varepsilon_1}{\delta(\varepsilon_1^2+\varepsilon_2+1)}, ~\frac{(\delta^2+\delta-1)\varepsilon_1\varepsilon_2} {\delta^3(1+\varepsilon_2)} \right\}
\ee
for any $\delta\in]\frac{\sqrt{5}-1}{2},+\infty[ $ to   ensure some convergence properties. Then    we get $\max\limits_{\delta\in]\frac{\sqrt{5}-1}{2},+\infty[}\kappa(\delta) =\kappa(\sqrt{3}-1)=\frac{1}{2}$, and use $\alpha\in]0,\kappa(\delta)[$ to control the step size. Our range of $\alpha$  is larger than that presented in \cite{9.,yang},  see  Lemma \ref{lem_kappa} for more explanations. Secondly, the region of $\delta$ is partitioned as
$$
\delta\in \left](\sqrt{5}-1)/2,+\infty\right[=\left](\sqrt{5}-1)/2,1 \right[~ \cup ~[1,+\infty[
$$
to explore  convergence of  the proposed method, and the  $\cO(1/n)$  ergodic convergence rate  is established. Finally, we obtain the upper bound of $\alpha$ by an extremely simple example, and improve numerical efficiency by introducing nonmonotonic step size $\lambda_n$ but $\frac{\lambda_n}{\lambda_{n-1}}\rightarrow1$.
In fact, the proposed nonmonotonic step size can break away from overdependence on the initial point, but it would have to be monotonic in the end for getting convergence.

The paper is organized as follows. In Section \ref{sec_preliminarries}, we provide some useful facts and notations. In Section \ref{sec_PEG}, we introduce our algorithm and explore the properties of the function $\kappa(\delta)$. A weak convergence theorem of our method is proved in Section \ref{sec_Convergence}. In Section \ref{sec_Rate}, we establish the ergodic convergence rate of the proposed algorithms, and we improve the algorithms in Section \ref{sec_improved} to avoid the adverse effects of the nonincreasing step size. In Section \ref{sec_fefurther}, we show by  an example  that any value of $\alpha\in ]\frac{2}{2\delta+1},+\infty[$ with $\delta\in]0,+\infty[$ does not guarantee convergence of our algorithm.  Numerical experiments on solving some problems tested in the   literatures  are provided and analyzed in Section \ref{sec_experiments}. We finally conclude our paper in Section \ref{sec_conclusion}.

\section{Preliminaries}
\label{sec_preliminarries}
In this section, we introduce some notations and facts on the well-known properties of the proximal operator, Opial condition and Young's inequality, which are used for  the sequel convergence analyses.

The proximal operator prox$_{\lambda g}:\cH\rightarrow \cH$ with prox$_{\lambda g}(x) = (I+\lambda \partial g)^{-1}(x), \lambda>0, x\in\cH$, is defined by
\ben
\mbox{prox}_{\lambda g}(x)
:= \argmin_{y\in\cH}\left\{ g(y)+\frac{1}{2\lambda}\|x-y\|^2\right\},\quad\forall\, x\in\cH, \lambda>0.
\een

Setting
\be\label{phi}
\Phi(x, y) := \langle F(x), y-x\rangle + g(y)-g(x),
\ee
it is clear that problem (\ref{primal_problem0}) is equivalent to finding $x^*\in\cH$ such that $\Phi(x^*, y)\geq0$ for all $y\in\cH$.

\begin{fact}\cite{Bauschke2011Convex}\label{fact_proj}
Let $g :\cH\rightarrow (-\infty, +\infty]$ be a convex function, $\lambda>0$ and $x\in \cH$. Then $p = \mbox{prox}_{\lambda g}(x)$ if and
only if
\ben
\langle p-x, y-p\rangle\geq \lambda [g(p)-g(y)],~~ \forall y\in\cH.
\een
\end{fact}

\begin{fact}\cite{Opial}\label{fact_weakly}
(Opial 1967) Let $\cS$ be a nonempty set of $\cH$ and $\{x_n\}_{k\in\bN}$ be a sequence in $\cH$ such that the following two conditions hold:\\
(1) for every $x^*\in \cS$, $\lim\limits_{n\rightarrow+\infty} \|x_n- x^*\|$ exists;\\
(2) every sequential weak cluster point of $\{x_n\}_{k\in\bN}$ is in $\cS$.\\
Then $\{x_n\}_{k\in\bN}$ converges weakly to a point in $\cS$.
\end{fact}

\begin{fact}\label{fact_ab}
Let $\{a_n\}$, $\{b_n\}$ be two nonnegative real sequences and $\exists N>0$ such that
\ben
a_{n+1} \leq a_n-b_n,~~\forall n>N.
\een
Then $\{a_n\}$ is convergent and $\lim\limits_{n\rightarrow \infty} b_n = 0$.
\end{fact}

\begin{fact}\label{fact_Yang}
(Young's inequality)  For all $a, b\geq0$ and $\varepsilon> 0$, we have
$$
ab\leq \frac{a^2}{2\varepsilon} + \frac{\varepsilon b^2}{2}.
$$
\end{fact}

The following identity (cosine rule) appears in many times and we will use it for simplicity of  convergence analyses. For all $x, y, z\in\cH$,
\be\label{id}
\langle x-y, x-z\rangle= \frac{1}{2}\|x-y\|^2 +\frac{1}{2}\|x-z\|^2- \frac{1}{2}\|y-z\|^2.
\ee

\section{Proximal Extrapolated Gradient Method with Prediction and Correction}
\label{sec_PEG}
In this section, we state our proximal extrapolated gradient method with prediction and correction (PEG), by using the step size function $\kappa(\delta)$ defined in (\ref{kappa}).

\vskip5mm
\hrule\vskip2mm
\begin{algo}
[PEG for solving (\ref{primal_problem0})]\label{algo1}
{~}\vskip 1pt {\rm
\begin{description}
\item[{\em Step 0.}] Take $\delta\in]\frac{\sqrt{5}-1}{2},+\infty[ $,  choose  $x_0\in \mathcal{H},$  $\lambda_0>0$, $\gamma\in(0,1)$, $\alpha \in ]0, \kappa(\delta)[$ and a bounded sequence $\{\zeta_n>0\}$. Set $y_0=x_0$, $x_{1}=\mbox{prox}_{\lambda_0 g}(x_0-\lambda_0 F(x_0))$ and $n=1$.
\item[{\em Step 1.}] \textbf{Prediction:}\\
1.a. Compute
\be
y_{n}&=&x_{n}+\delta (x_{n}-x_{n-1}),\\
\lambda_{n}&=&
         \min~\left\{\lambda_{n-1}, ~~\frac{\alpha\|y_{n}-y_{n-1}\|} {\|F(y_{n})-F(y_{n-1})\|} \right\}.\label{lambda_tilde}
\ee
1.b. Compute \ \ \ \
\ben
 x_{n+1}=\mbox{prox}_{\lambda_n g}(x_n-\lambda_n F(y_n)),
\een
if $x_{n+1}=x_n=y_n$, then stop: $x_{n+1} $ is a solution.
\item[{\em Step 2.}] \blue{\textbf{Correction when $\delta<1$:}\par
Check
\ben
\|x_{n+1}-x_n\|\leq \zeta_n,
\een
if not hold, set $\lambda_n\leftarrow \gamma \lambda_n$ and return to Step 1.b.}
\item[{\em Step 3.}] Set $n\leftarrow n + 1$ and return to Step 1.\\
  \end{description}
}
\end{algo}
\vskip1mm\hrule\vskip5mm

The aim of Correction step is to bound $\{\|x_{n}-x_{n-1}\|\}$ by the given sequence $\{\zeta_n\}$ when $\delta<1$, as convergence analysis requires $\|x_{n+1}-x_{n}\|<+\infty$. In practice, we don't need to give the sequence $\{\zeta_n\}$, but generate adaptively by
\be\label{zeta}
\zeta_n=\max\{\zeta_{\min},~~\min\{\mu \|x_{n}-x_{n-1}\|,~~\nu\|x_{1}-x_{0}\|\}\},
\ee
for given $1<\mu\leq\nu$ and small $\zeta_{\min}$ (e.g., $\zeta_{\min}=10^{-6}$), then $\zeta_n\leq \nu\|x_{1}-x_{0}\|$ for all $n\geq 1$ and $\zeta_n\geq \zeta_{\min}$. Moreover, we observe $\|x_{n+1}-x_n\|\leq\mu \|x_{n}-x_{n-1}\|$ for bounding more tightly due to $\|x_{n+1}-x_n\|\rightarrow 0$.

For a convex function $f$, if $F=\nabla f$ we observe $\|x_{n+1}-x_{n}\|<+\infty$, see (\ref{bounded}), so Correction step is not necessary. However for other cases, one needs to apply linesearch to ensure $\|x_{n+1}-x_{n}\|<+\infty$. Interestingly, for all the tested problems shown in Section \ref{sec_experiments}, the linesearch in Correction step does not start to arrive termination conditions, when using (\ref{zeta}) with $\mu=\nu=10$. Namely, the predicted step is good enough for obtaining a convergent sequence for the tested problems, though the convergence without prediction is unknown in general.

The following lemma shows that the correction procedure described in Algorithm \ref{algo1} is well-defined.

\begin{lem}\label{lem_cor}
The correction procedure always terminates. i.e., $\{\lambda_n\}$ is well defined when $\delta\in]\frac{\sqrt{5}-1}{2}, 1[$.
\end{lem}
\proof
Denote
\ben
A:=\partial g ~~\mbox{and}~~x_{n+1}(\lambda) := \mbox{prox}_{\lambda g}(x_n-\lambda F(y_n)).
\een
From \cite[Theorem 23.47]{Bauschke2011Convex}, we have that $\mbox{prox}_{\lambda g}[x_{n+1}(0)]\rightarrow P_{\overline{\dom A}}[x_{n+1}(0)]$ as $\lambda\rightarrow 0$ ($\overline{\dom A}$ denotes the
closures of $\dom A$), which together with the nonexpansivity of $\mbox{prox}_{\lambda g}$ yields
\ben
&&\|x_{n+1}(\lambda)-P_{\overline{\dom A}}[x_{n+1}(0)]\|\\
&\leq& \|x_{n+1}(\lambda)-\mbox{prox}_{\lambda g}[x_{n+1}(0)]\|+\|\mbox{prox}_{\lambda g}[x_{n+1}(0)] -P_{\overline{\dom A}}[x_{n+1}(0)]\|\\
&\leq&\lambda\|F(y_n)\|+\|\mbox{prox}_{\lambda g}[x_{n+1}(0)] -P_{\overline{\dom A}}[x_{n+1}(0)]\|.
\een
By taking the limit as $\lambda\rightarrow 0$, we deduce that $x_{n+1}(\lambda)\rightarrow P_{\overline{\dom A}}[x_{n+1}(0)]$. Notice that $x_{n+1}(0)=x_n$, we observe $P_{\overline{\dom A}}[x_{n+1}(0)]=x_n$.

By a contradiction, suppose that the correction procedure in Algorithm \ref{algo1} fails to terminate at the $n$-th iteration.  Then, for all $\lambda= \gamma^i \lambda_n$ with $i=0,1,\cdots $, we have $\|x_{n+1}(\lambda)-x_n\| >\zeta_n$.
Since $\gamma^i\rightarrow 0$ as $i\rightarrow\infty$, so $\lambda\rightarrow 0$, this gives a contradiction $0\geq\zeta_n$, which completes the proof.
\endproof

\begin{rem}\label{rem_lower_bound}
Note that the sequence $\{\lambda_n\}$ is monotonically decreasing. Since  $F$ is a $L$-Lipschitz continuous mapping ($L>0$),  we have
\ben
\frac{\alpha\|y_n-y_{n-1}\|}{\|F(y_n)-F(y_{n-1})\|}\geq\frac{\alpha\|y_n-y_{n-1}\|}{L\|y_n-y_{n-1}\|}=\frac{\alpha}{L} \een
for $F(y_n)\neq F(y_{n-1})$. Thus the predicted step sequence $ \{\lambda_n\}_{n\in\bN} $ has a lower bound $\tau:=\min\{{\frac{\alpha}{L },\lambda_0 }\}$,  then when $\delta\geq1$ its limit exists and $\lim\limits_{n\rightarrow\infty}\lambda_n\geq\tau>0$. If $\delta<1$, $\{\lambda_n\}$ is well defined from Lemma \ref{lem_cor}, and has a lower bound $\tau:=\min\{{\frac{\gamma^{i_0}\alpha}{L },\lambda_0 }\}$ for some $i_0\geq0$, which implies $\lim\limits_{n\rightarrow\infty}\lambda_n>0$ as well.
\end{rem}

Below, we derive the analytical expression of $\kappa(\delta)$.

\begin{lem}\label{lem_kappa}
For the function $\kappa(\delta)$ defined in (\ref{kappa}), we have $\kappa(\delta)=\frac{\sqrt{a+1}}{\delta(a+1+\sqrt{a+1})}$ with $a=\frac{\delta^2}{\delta^2+\delta-1}$ for $\delta\in]\frac{\sqrt{5}-1}{2},+\infty[$.
\end{lem}
\proof
Fix  $\delta\in]\frac{\sqrt{5}-1}{2},+\infty[$, then $\delta^2+\delta-1>0$.  Noting that the structure of (\ref{kappa}) and $\kappa(\delta)$ is a maximum value, so $\frac{\varepsilon_1}{\delta(\varepsilon_1^2+\varepsilon_2+1)}= ~\frac{(\delta^2+\delta-1)\varepsilon_1\varepsilon_2} {\delta^3(1+\varepsilon_2)}$, which together with $a=\frac{\delta^2}{\delta^2+\delta-1}$ and $\varepsilon_1=\sqrt{a+1}$ shows
\be\label{k_2}
\kappa(\delta)=\max\limits_{\varepsilon_2>0} \frac{\sqrt{a\varepsilon_2+(a-1) \varepsilon_2^2-\varepsilon_2^3}} {\delta(a+a\varepsilon_2)}.
\ee
By the first-order optimality condition of the optimization problem (\ref{k_2}), we have $\varepsilon_2=\sqrt{a+1}-1$. Substituting it into (\ref{k_2}), the result can be deduced.
\endproof

\begin{figure}[htpb]
\centering
\includegraphics[width=0.65\textwidth]{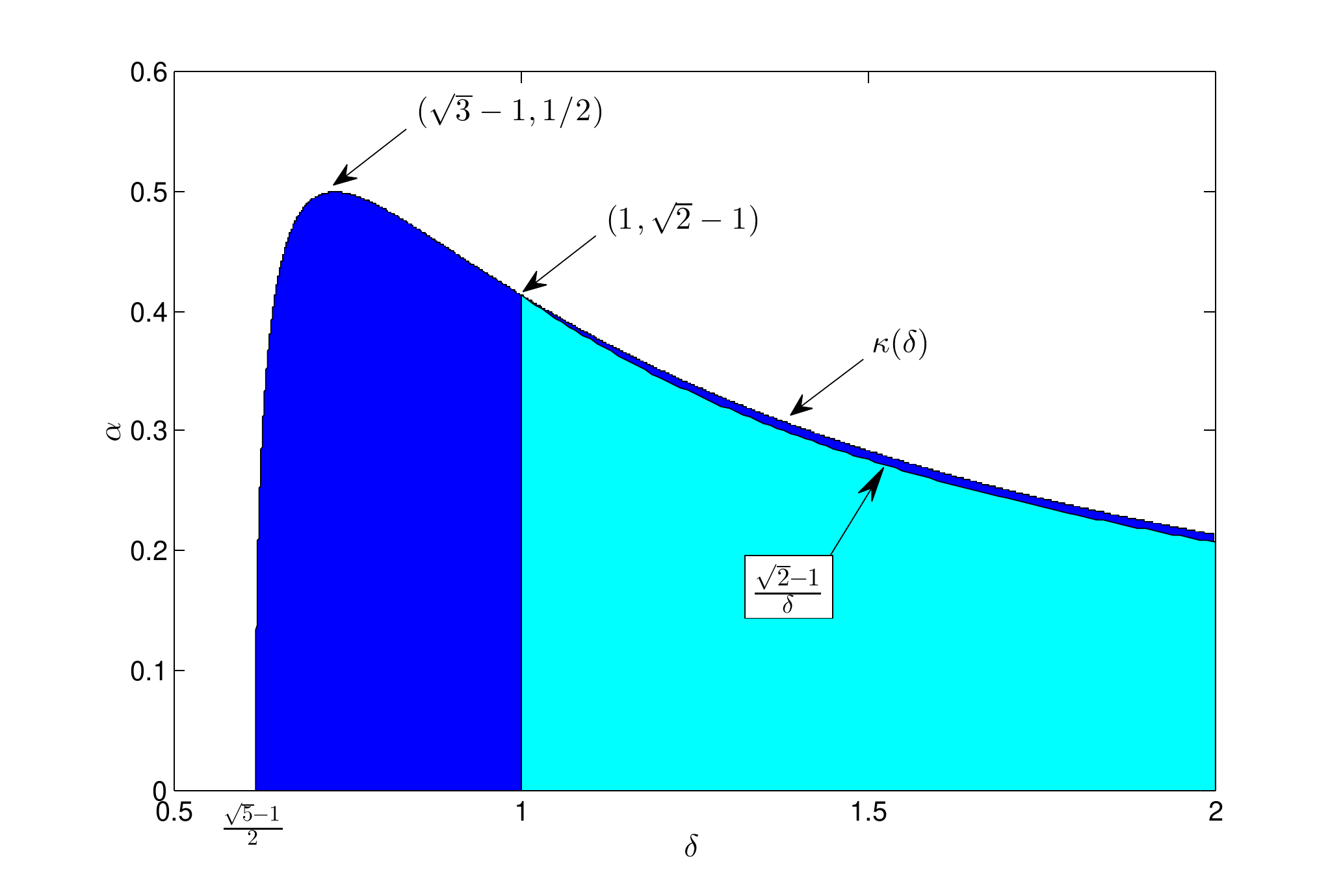}
\caption{The improved region (blue) of the parameters $\alpha$ and $\delta$.}
\label{Fig 1} 
\end{figure}

By Lemma \ref{lem_kappa} and Fig. \ref{Fig 1}, the maximum value of $\kappa(\delta)$ is $\frac{1}{2}$ when $\delta\in]\frac{\sqrt{5}-1}{2},+\infty[$, and in fact,
$$
\max\limits_{\delta\in]\frac{\sqrt{5}-1}{2},+\infty[}\kappa(\delta)
=\kappa(\sqrt{3}-1)=\frac{1}{2}.
$$
In this case, we have $a=2$, $\varepsilon_1=\sqrt{3}$ and $\varepsilon_2=\sqrt{3}-1$.

\begin{rem}\label{rem_extent_range}
It can be noticed that the method proposed in \cite{yang} is a special case of Algorithm \ref{algo1}, when $g(x)=l_{C}(x)$ and $\delta\in]1,+\infty[ $, but $\kappa(\delta) > \frac{\sqrt{2}-1}{\delta }$ from Lemma \ref{lem_kappa}. Namely, we extend the range of $\delta$ and then improve the upper bound of $\alpha$ when the operator is the gradient of a convex function or using linesearch, see Fig. \ref{Fig 1}, which causes larger step size $\lambda_n$ that will be more efficient for numerical experiments.
\end{rem}

\subsection{Convergence Analysis}
\label{sec_Convergence}
This section devotes to studying    convergence properties of Algorithm \ref{algo1}. For $\delta\in[1,+\infty[ $, its   convergence and convergence rate can be obtained by combining the methods in \cite{yang,Proximal-extrapolated} with  the basic theory of limit. However, it is a completely different situation for  $\delta\in]\frac{\sqrt{5}-1}{2},1[$, since the desired properties (such as monotonicity and nonnegativity) are no longer valid in the case of $\delta\in]\frac{\sqrt{5}-1}{2},1[$  although we can adopte a larger value of $\alpha$.

We next give a basic lemma about the iterations generated by  Algorithm \ref{algo1}     for any $\delta\in]\frac{\sqrt{5}-1}{2},+\infty[$, which play a crucial role in  proving    the main convergence results.

\begin{lem}\label{lem0}
Let $\{x_n\}$ and $\{y_n\}$ be two sequences generated by Algorithm \ref{algo1}. For any  $x\in \cH$, we have
\be\label{result_lem1}
\|x_{n+1}-x\|^2&\leq&\|x_n-x\|^2-\frac{\lambda_{n}}{\delta \lambda_{n-1}}[\|y_n-x_n\|^2+\|x_{n+1}-y_n\|^2]+
\left(\frac{\lambda_n}{\delta \lambda_{n-1}}-1 \right)\|x_{n+1}-x_n\|^2   \nonumber\\
&&+ 2\alpha\|y_n-y_{n-1}\|\|x_{n+1}-y_n\|-2\lambda_n[(1+ \delta)g(x_n)- \delta g(x_{n-1})-g(x)]\nonumber\\
&&-2\lambda_n \langle F(y_n), y_n-x\rangle.
\ee
\end{lem}
\proof
Followed by $x_{n+1}=\mbox{prox}_{\lambda_n g}(x_n-\lambda_n F(y_n))$ and Fact \ref{fact_proj}, we have
\be\label{temp_01}
\langle x_{n+1}-x_{n} + \lambda_{n}F(y_{n}), x-x_{n+1}\rangle\geq \lambda_{n}(g(x_{n+1})-g(x)),~~ \forall x\in\cH,
\ee
which shows
\ben
\langle x_{n}-x_{n-1} + \lambda_{n-1}F(y_{n-1}), x-x_{n}\rangle\geq \lambda_{n-1}(g(x_{n})-g(x)),~~ \forall x\in\cH.
\een
Substituting $x:= x_{n+1}$ and $x:= x_{n-1}$ into the above inequality respectively, we obtain
\be
\langle x_{n}-x_{n-1} + \lambda_{n-1}F(y_{n-1}), x_{n+1}-x_{n}\rangle\geq \lambda_{n-1}(g(x_{n})-g(x_{n+1})),\label{tem1}\\
\langle x_{n}-x_{n-1} + \lambda_{n-1}F(y_{n-1}), x_{n-1}-x_{n}\rangle\geq \lambda_{n-1}(g(x_{n})-g(x_{n-1})).\label{tem2}
\ee
Multiplying (\ref{tem2}) by $\delta$ and then adding it to (\ref{tem1}), which by $y_{n}=x_{n}+\delta (x_{n}-x_{n-1})$  yields
\be\label{tem3}
\langle x_{n}-x_{n-1} + \lambda_{n-1}F(y_{n-1}), x_{n+1}-y_{n}\rangle\geq \lambda_{n-1}[(1+\delta)g(x_{n})-g(x_{n+1})-\delta g(x_{n-1})].
\ee
Multiplying (\ref{tem3}) by $\frac{\lambda_{n}}{\lambda_{n-1}}$ and using $y_{n}=x_{n}+\delta (x_{n}-x_{n-1})$ again, we get
\be\label{temp_02}
\left\langle \frac{\lambda_n(y_{n}-x_{n})}{\delta\lambda_{n-1}} + \lambda_{n}F(y_{n-1}), x_{n+1}-y_{n}\right\rangle\geq \lambda_{n}[(1+\delta)g(x_{n})-g(x_{n+1})-\delta g(x_{n-1})].
\ee
Finally, adding (\ref{temp_01}) to (\ref{temp_02}) gives us
\ben
&&\langle x_{n+1}-x_{n}, x-x_{n+1}\rangle+\left\langle \frac{\lambda_n(y_{n}-x_{n})}{\delta\lambda_{n-1}}, x_{n+1}- y_n\right\rangle + \lambda_n\langle F(y_n)-F(y_{n-1}), y_n-x_{n+1}\rangle\\
&\geq&\lambda_{n}[(1+\delta)g(x_{n})-g(x)-\delta g(x_{n-1})]+\lambda_n \langle F(y_n), y_n-x\rangle.
\een
Then, using (\ref{id}), the updating of $\lambda_n$ and Cauchy-Schwarz inequality, we obtain
\ben
\|x_{n+1}-x\|^2&\leq&\|x_n-x\|^2-\frac{\lambda_{n}}{\delta \lambda_{n-1}}[\|y_n-x_n\|^2+\|x_{n+1}-y_n\|^2]+
\left(\frac{\lambda_n}{\delta \lambda_{n-1}}-1 \right)\|x_{n+1}-x_n\|^2   \\
&&+ 2\lambda_n\|F(y_n)-F(y_{n-1})\|\|x_{n+1}-y_n\|-2\lambda_n[(1+ \delta)g(x_n)- \delta g(x_{n-1})-g(x)]\\
&&-2\lambda_n \langle F(y_n), y_n-x\rangle\\
&\leq&\|x_n-x\|^2-\frac{\lambda_{n}}{\delta \lambda_{n-1}}[\|y_n-x_n\|^2+\|x_{n+1}-y_n\|^2]+
\left(\frac{\lambda_n}{\delta \lambda_{n-1}}-1 \right)\|x_{n+1}-x_n\|^2   \\
&&+ 2\alpha\|y_n-y_{n-1}\|\|x_{n+1}-y_n\|-2\lambda_n[(1+ \delta)g(x_n)- \delta g(x_{n-1})-g(x)]\\
&&-2\lambda_n \langle F(y_n), y_n-x\rangle.
\een
The proof is completed.
\endproof

\begin{lem}\label{lem1}
 Let $\{x_n\}$, $\{y_n\}$ be two sequences generated by Algorithm \ref{algo1} and $\bar{x}\in \cS$ (the solution set of problem (\ref{primal_problem0})). Then, for any $\varepsilon_1,\varepsilon_2>0$, we have
\ben
\|x_{n+1}-\bar{x}\|^2+2\lambda_n(1 +\delta)	\Phi(\bar{x}, x_n)&\leq&\|x_n-\bar{x}\|^2+ 2\lambda_{n-1}(1 +\delta)	\Phi(\bar{x}, x_{n-1})\\
&&+\left[\frac{1}{\varepsilon_1}\left(1+\frac{1}{\varepsilon_2}\right)\alpha-\frac{\lambda_n}{\delta \lambda_{n-1}} \right]\|x_n-y_n\|^2\\
&&+
\left(\frac{\lambda_n}{\delta \lambda_{n-1}}-1 \right)\|x_{n+1}-x_n\|^2   \\
&&+\left(\varepsilon_1 \alpha-\frac{\lambda_n}{\delta \lambda_{n-1} }\right)\|x_{n+1}-y_n\|^2+ \frac{1+\varepsilon_2}{\varepsilon_1}\alpha\|x_n-y_{n-1}\|^2,
\een
where $\Phi$ is defined as in (\ref{phi}).
\end{lem}
\proof
Using Fact \ref{fact_Yang}, for any $\varepsilon_1>0$  we have
\ben
2\alpha\|y_n-y_{n-1}\|\|y_n-x_{n+1}\|
 &\leq&\alpha(\frac{1}{\varepsilon_1 }\|y_n-y_{n-1}\|^2+\varepsilon_1\|x_{n+1}-y_n\|^2).
\een
Meanwhile, for any $\varepsilon_2>0$  we deduce
\ben
\|y_n-y_{n-1}\|^2&=&\|y_n-x_n\|^2+\|x_n-y_{n-1}\|^2+2\langle y_n-x_n,x_n-y_{n-1}\rangle   \\
 &\leq&\|y_n-x_n\|^2+\|x_n-y_{n-1}\|^2+2\|y_n-x_n\|\|x_n-y_{n-1}\|  \\
 &\leq&(1+\frac{1}{\varepsilon_2})\|y_n-x_n\|^2+(1+\varepsilon_2) \|x_n-y_{n-1}\|^2.
\een
Combining the above inequalities we have
\be
&&2\alpha\|y_n-y_{n-1}\|\|y_n-x_{n+1}\|  \nonumber\\
&\leq&\alpha\left[\frac{1}{\varepsilon_1}\left(1+\frac{1}{\varepsilon_2}\right)\|y_n-x_n\|^2 +\frac{1+\varepsilon_2}{\varepsilon_1} \|x_n-y_{n-1}\|^2+\varepsilon_1\|x_{n+1}-y_n\|^2\right].\label{ineq2}
\ee
In addition, the monotonicity of $F$ implies for any $x\in\cH$
\be\label{ineq1}
\lambda_n\langle F(y_n),y_n-x\rangle &\geq& \lambda_n\langle F(x),y_n-x\rangle \nonumber\\
&=&\lambda_n[(1+\delta)\langle F(x),x_n-x\rangle-\delta \langle F(x),x_{n-1}-x\rangle].
\ee

Substituting (\ref{ineq2}) and (\ref{ineq1}) into (\ref{result_lem1}), we deduce by the aids of $\Phi(x,y)$ in  (\ref{phi}) that
\be\label{rate}
\|x_{n+1}-x\|^2&\leq&\|x_n-x\|^2 -\|x_{n+1}-x_n\|^2+\frac{1}{\varepsilon_1}\left(1+\frac{1}{\varepsilon_2}\right) \alpha\|y_n-x_n\|^2\nonumber\\
&&+\frac{1+\varepsilon_2}{\varepsilon_1}\alpha\|x_n-y_{n-1}\|^2+\varepsilon_1\alpha\|x_{n+1}-y_n\|^2\nonumber\\
 &&+ \frac{\lambda_n}{\delta \lambda_{n-1} }(\|x_{n+1}-x_n\|^2-\|x_n-y_n\|^2-\|x_{n+1}-y_n\|^2)\nonumber\\
 &&- 2\lambda_n[(1 +\delta)	\Phi(x, x_n)-\delta	\Phi(x, x_{n-1})].
\ee
Since $\delta>0$ and $ \{\lambda_n\}_{n\in\bN}$ is a monotone decreasing sequence, we have $\lambda_n \delta \leq \lambda_{n-1} \delta\leq (1 + \delta)\lambda_{n-1}$. Note that $\Phi(\bar{x}, x_{n-1})\geq0$ for any $\bar{x}\in \cS$, then
\ben
\|x_{n+1}-\bar{x}\|^2&\leq&\|x_n-\bar{x}\|^2+\left[\frac{1}{\varepsilon_1}\left(1+\frac{1}{\varepsilon_2}\right)\alpha-\frac{\lambda_n}{\delta \lambda_{n-1}} \right]\|x_n-y_n\|^2\\
&&+
\left(\frac{\lambda_n}{\delta \lambda_{n-1}}-1 \right)\|x_{n+1}-x_n\|^2   \\
&&+\left(\varepsilon_1 \alpha-\frac{\lambda_n}{\delta \lambda_{n-1} }\right)\|x_{n+1}-y_n\|^2+ \frac{1+\varepsilon_2}{\varepsilon_1}\alpha\|x_n-y_{n-1}\|^2\\
&&- 2\lambda_n(1 +\delta)	\Phi(\bar{x}, x_n)+ 2\lambda_{n-1}(1 +\delta)	\Phi(\bar{x}, x_{n-1}).
\een
This completes the proof.\qed

By Lemma \ref{lem1} and some transpositions, we have the following results directly.

\begin{lem}\label{lem4}
Let $\{x_n\}$, $\{y_n\}$ be two sequences generated by Algorithm \ref{algo1} and $\bar{x}\in \cS$. Then, for any $\varepsilon_1,\varepsilon_2>0$,
we have
\be\label{eq_thm}
a_{n+1}&\leq& a_n-b_n,
\ee
where
\be\label{a_n1}
\left\{
\ba{rcl}
a_n&=&\|x_n-\bar{x}\|^2+ 2\lambda_{n-1}(1 +\delta)	\Phi(\bar{x}, x_{n-1})+\frac{1+\varepsilon_2}{\varepsilon_1}\alpha\|x_n-y_{n-1}\|^2\\
&&+\left(1- \frac{\lambda_{n-1}}{\delta \lambda_{n-2}} \right)\|x_n-x_{n-1}\|^2, ~~n\geq2,   \\
b_n&=&\left[\frac{\lambda_n}{\delta \lambda_{n-1} }-\left(\varepsilon_1+\frac{1+\varepsilon_2}{\varepsilon_1}\right) \alpha\right]\|x_{n+1}-y_n\|^2\\
&&+\left[\frac{\lambda_n}{\delta \lambda_{n-1}}-\frac{1}{\varepsilon_1}\left(1+\frac{1}{\varepsilon_2}\right)\alpha +\frac{1}{\delta^2}\left(1- \frac{\lambda_{n-1}}{\delta \lambda_{n-2}} \right)\right]\|x_n-y_n\|^2,
\ea
\right.
\ee
or
\be\label{a_n2}
\left\{\ba{rcl}
a_n&=&\|x_n-\bar{x}\|^2+ 2\lambda_{n-1}(1 +\delta)	\Phi(\bar{x}, x_{n-1})+\frac{1+\varepsilon_2}{\varepsilon_1}\alpha\|x_n-y_{n-1}\|^2, ~~n\geq2,   \\
b_n&=&\left[\frac{\lambda_n}{\delta \lambda_{n-1} }-\left(\varepsilon_1+\frac{1+\varepsilon_2}{\varepsilon_1}\right) \alpha\right]\|x_{n+1}-y_n\|^2\\
&&+\left[\frac{\lambda_n}{\delta \lambda_{n-1}}-\frac{1}{\varepsilon_1}\left(1+\frac{1}{\varepsilon_2}\right)\alpha \right]\|x_n-y_n\|^2
+\left(1- \frac{\lambda_n}{\delta \lambda_{n-1}} \right)\|x_{n+1}-x_n\|^2.
\ea\right.
\ee
\end{lem}

Because the sequence $\{\lambda_n\}_{n\in\bN}$ is monotonically decreasing, we have $1-\frac{\lambda_n}{\delta \lambda_{n-1}}\geq0$ for any $\delta\geq1$. But for $\delta\in]\frac{\sqrt{5}-1}{2},1[$, we have $\lim\limits_{n\rightarrow+\infty}\left(1-\frac{\lambda_n}{\delta \lambda_{n-1}} \right)=1-\frac{1}{\delta}<0$. So,   convergence of Algorithm \ref{algo1} with $\delta<1$ is different from that with $\delta\geq1$, and hence cannot be established by the similar methods as in \cite{yang,Proximal-extrapolated}.

Notice that $a_n\geq0$ in (\ref{a_n1}) when $\delta\geq1$, we take (\ref{a_n1}) to study the convergence of Algorithm \ref{algo1} with $\delta\geq1$. Consequently, a larger upper bound $\kappa(\delta)$ of $\alpha$ is obtained than that in \cite{yang}. While for the case of  $\delta<1$, we take (\ref{a_n2}) as $a_n\geq0$ for all $n\geq1$, and further investigate the properties of $b_n$ to ensure convergence of Algorithm \ref{algo1}.

Below we state and prove our main convergence result of Algorithm \ref{algo1} for above two different regions: $\delta\in]\frac{\sqrt{5}-1}{2},1[$ and $\delta\in[1,+\infty[$.

\begin{thm}\label{thm21}
Let $\{x_n\}$ be the sequence generated by Algorithm \ref{algo1} with $\delta\in[1,+\infty[$. Then, $\{x_n\}$ converges weakly to a solution of problem (\ref{primal_problem0}).
\end{thm}
\proof
From Remark \ref{rem_lower_bound}, we have $\lim\limits_{n\rightarrow\infty}\lambda_n=\lambda>0$. Then for any $\delta\in[1,+\infty[$ and  $\alpha<\kappa(\delta)$, we have
\ben
&&\lim\limits_{n\rightarrow\infty}\left[\frac{\lambda_n}{\delta \lambda_{n-1}}-\left(\varepsilon_1+\frac{1+\varepsilon_2}{\varepsilon_1} \right)\alpha\right]
=\frac{1}{\delta }-\left(\frac{\varepsilon_1^2+\varepsilon_2+1}{\varepsilon_1}\right) \alpha>0,\\
&&\lim\limits_{n\rightarrow\infty}\left[\frac{\lambda_{n}}{\delta \lambda_{n-1} }-\frac{1}{\varepsilon_1}\left(1+\frac{1}{\varepsilon_2}\right)\alpha +\frac{1}{\delta^2}\left(1- \frac{\lambda_{n-1}}{\delta \lambda_{n-2}} \right)\right]\\
&=&\frac{1}{\delta }-\frac{1}{\varepsilon_1}\left(1+\frac{1}{\varepsilon_2}\right)\alpha +\frac{1}{\delta^2}\left(1- \frac{1}{\delta} \right) >0.
\een
Thus, there exists an integer $N>2,$ such that for any $n>N$,
\ben
\left.
\ba{r}
\frac{\lambda_n}{\delta \lambda_{n-1} }-\left(\frac{\varepsilon_1^2+\varepsilon_2+1}{\varepsilon_1}\right) \alpha>0,\\
\frac{\lambda_{n}}{\delta \lambda_{n-1} }-\frac{1}{\varepsilon_1}\left(1+\frac{1}{\varepsilon_2}\right)\alpha+\frac{1}{\delta^2}\left(1- \frac{\lambda_{n-1}}{\delta \lambda_{n-2}} \right) >0,
\ea
\right\}
\een
which implies that $b_n\geq0$ in (\ref{a_n1}) when $n>N$. Recall $1-\frac{\lambda_{n-1}}{\delta\lambda_{n-2}}\geq0$ for any $\delta\geq1$, we deduce $a_n\geq0$ in (\ref{a_n1}).  Hence, by Lemma \ref{lem4} and Fact \ref{fact_ab},  $\{a_n\}_{n\in\bN}$ is convergent and
$\lim\limits_{n\rightarrow\infty} b_n = 0$. This means that $\{
\|x_n-\bar{x}\|^2\}$ is bounded and so does $\{x_n\}_{n\in\bN}$. Also, we have   $\lim\limits_{n\rightarrow\infty}\|x_{n+1}-y_n\|=0$ and $\lim\limits_{n\rightarrow\infty}\|x_n-y_n\|=0.$ By $\|x_{n+1}-x_n\|=\frac{1}{\delta}\|x_{n+1}-y_{n+1}\|$, we also have that
$\lim\limits_{n\rightarrow\infty}\|x_{n+1}-x_n\|=0$ and $\{
y_n\}_{n\in\bN}$ is bounded.

In what follows, we   prove the sequence $\{x_{n}\}$ converges weakly to a solution of problem (\ref{primal_problem0}).
For any   cluster $x^*\in \cH$ of $\{x_{n}\}$, there exists a subsequence $\{x_{n_{k}}\}$ that converges weakly to $x^*$, namely $x_{n_{k}}\rightharpoonup x^*$. It is obvious that $\{y_{n_{k}}\}$ also converges weakly to $x^*$.  Next we verify that $x^*\in \cS$. Applying  Fact \ref{fact_proj}, we deduce
\be\label{ineq3}
\left\langle \frac{x_{n_k+1}-x_{n_k}}{\lambda_{n_k}}+ F(y_{n_k}), x-x_{n_k+1}\right\rangle\geq g(x_{n_k+1})-g(x),~~ \forall x\in\cH.
\ee
Letting $k\rightarrow\infty$ in (\ref{ineq3}) and using the facts $\lim\limits_{k\rightarrow\infty}\|x_{n_k+1}-x_{n_k}\| =0$, $g(x)$ is lower semicontinuous and $\lim\limits_{n\rightarrow\infty}\lambda_n=\lambda>0,$ we obtain
\ben
\langle F(x^*), x-x^*\rangle \geq \lim \inf\limits_{k\rightarrow\infty}g(x_{n_k+1})-g(x)\geq g(x^*)-g(x),~~ \forall x\in\cH,
\een
which confirms $x^*\in \cS$.

Finally, we prove that $x_{n}\rightharpoonup x^*$. We take $\bar{x}=x^*$ in the definition (\ref{a_n1}) of $a_n$ and label as $a_n^*$. Notice that $\{\lambda_n\}$ is bounded and $\Phi(x^*, \cdot)$ is continuous from   (\textbf{A3}), we observe
\ben
&&\lim\limits_{n\rightarrow\infty}a_{n}^*
=\lim\limits_{k\rightarrow\infty}a_{n_{k}+1}^* \\ &=&\lim\limits_{k\rightarrow\infty}\left(\ba{l}\|x_{n_k+1}-x^*\|^2+ 2\lambda_{n_k}(1 +\delta)	\Phi(x^*, x_{n_k})+\frac{1+\varepsilon_2}{\varepsilon_1} \alpha\|x_{n_k+1}-y_{n_k}\|^2\\+\left(1- \frac{\lambda_{n_k}}{\delta \lambda_{n_k-1}} \right)\|x_{n_k+1}-x_{n_k}\|^2\ea\right)\\
&=& 0.
\een
Therefore, $\lim\limits_{k\rightarrow\infty}\|x_{n}-x^*\|=0$, which  by Fact \ref{fact_weakly} shows  $x_{n}\rightharpoonup x^*$.
\endproof

Now, we focus on   convergence analysis of Algorithm \ref{algo1} with $\delta\in]\frac{\sqrt{5}-1}{2},1[$ and use (\ref{a_n2}). For this case, we can not establish the nonnegativity of $\{b_n\}$ and the monotonic decreasing of $\{a_n\}$  because   $\frac{1}{\delta}-1>0$. Consequently,   convergence of $\{a_n\}_{n\in\bN}$ can not be obtained from (\ref{eq_thm}). We thus need to further investigate the sequence $\{a_n\}$ for getting a clear convergence, by using the boundedness of $\{\|x_{n}-x_{n-1}\|\}$ from Correction step.

First, we show that $\|x_{n+1}-x_{n}\|<+\infty$ when $\delta\in]\frac{\sqrt{5}-1}{2},1[$ and the operator $F$ is the gradient of a convex function $f:\cH\rightarrow\bR$, i.e., $F=\nabla f$. From (\ref{rate}) with $x=x_n$ and $\Phi(x_n, x_n)=0$, we deduce
\ben
&&\left[2-\frac{\lambda_n}{\delta \lambda_{n-1} }\right]\|x_{n+1}-x_n\|^2+\frac{1+\varepsilon_2}{\varepsilon_1}\alpha\|x_{n+1}-y_{n}\|^2\\
&\leq&\left[\frac{1}{\varepsilon_1}\left(1+\frac{1}{\varepsilon_2}\right) \alpha-\frac{\lambda_n}{\delta \lambda_{n-1} }\right]\|y_n-x_n\|^2+\frac{1+\varepsilon_2}{\varepsilon_1}\alpha\|x_n-y_{n-1}\|^2\nonumber\\
&& +\left[\left(\varepsilon_1+\frac{1+\varepsilon_2}{\varepsilon_1}\right)\alpha-\frac{\lambda_n}{\delta \lambda_{n-1} }\right]\|x_{n+1}-y_n\|^2+2\lambda_n\delta	\Phi(x_n, x_{n-1}).
\een
Using $F=\nabla f$ and the convexity of $f$ yields
\ben
\Phi(x_n, x_{n-1})&=&\langle \nabla f(x_n), x_{n-1}-x_n\rangle + g(x_{n-1})-g(x_n)\\
&\leq & \phi(x_{n-1})-\phi(x_{n}),
\een
where $\phi=f+g$. This together with $\lim\limits_{n\rightarrow\infty}\lambda_n=\lambda>0$, $\alpha<\kappa(\delta)$ and $\lambda_n\leq\lambda_{n-1}$ gives us
\ben
&&\left(2-\frac{1}{\delta}\right)\|x_{n+1}-x_n\|^2+\frac{1+\varepsilon_2}{\varepsilon_1}\alpha\|x_{n+1}-y_{n}\|^2+ 2\lambda_n\delta(\phi(x_{n})-\phi(\bar{x}))\\
&\leq&\frac{1+\varepsilon_2}{\varepsilon_1}\alpha\|x_n-y_{n-1}\|^2+ 2\lambda_{n-1}\delta(\phi(x_{n-1})-\phi(\bar{x})), ,
\een
where $\bar{x}\in\cS$, which implies from $\phi(x_{n})-\phi(\bar{x})\geq0$ that
\be\label{bounded}
\left(2-\frac{1}{\delta}\right)\|x_{n+1}-x_n\|^2\leq\frac{1+\varepsilon_2}{\varepsilon_1}\alpha \|x_{N+1}-y_{N}\|^2+ 2\lambda_{N}\delta(\phi(x_{N})-\phi(\bar{x}))<+\infty,~~\forall n>N.
\ee
That is to say, Correction step is not necessary when $F=\nabla f$, for a convex function $f$.

\begin{thm}\label{thm2}
Let $\{x_n\}$ be the sequence generated by Algorithm \ref{algo1} with $\delta\in]\frac{\sqrt{5}-1}{2},1[$. Then,  $\{x_n\}$ converges weakly to a solution of problem (\ref{primal_problem0}).
\end{thm}
\proof
Firstly, $\delta\in]\frac{\sqrt{5}-1}{2},1[$ gives $\frac{1}{\delta }>\frac{\delta^2 +\delta -1}{\delta^3}>0$. Note that $\lim\limits_{n\rightarrow\infty}\lambda_n=\lambda>0$, by taking the limit and from $\alpha<\kappa(\delta)$, we have
\be\label{ineq_6}
\left.
\ba{r}
\lim\limits_{n\rightarrow\infty}\left[\left(\varepsilon_1+\frac{1+\varepsilon_2}{\varepsilon_1} \right)\alpha-\frac{\lambda_n}{\delta \lambda_{n-1} }\right]
=\left(\frac{\varepsilon_1^2+\varepsilon_2+1}{\varepsilon_1}\right) \alpha-\frac{1}{\delta }<0,\\
\lim\limits_{n\rightarrow\infty}\left[\frac{1}{\varepsilon_1}\left(1+\frac{1}{\varepsilon_2}\right)\alpha -\frac{\lambda_{n}}{\delta \lambda_{n-1} }\right]=\frac{1}{\varepsilon_1}\left(1+\frac{1}{\varepsilon_2}\right)\alpha -\frac{1}{\delta }<0,\\
\lim\limits_{n\rightarrow\infty}\left[\frac{1}{\varepsilon_1} \left(1+\frac{1}{\varepsilon_2}\right)\alpha-\frac{\lambda_n}{\delta \lambda_{n-1} }+\frac{1}{\delta^2}\left(\frac{\lambda_{n-1}}{\delta \lambda_{n-2}}-1 \right)\right]
=\frac{1}{\varepsilon_1}\left(1+\frac{1}{\varepsilon_2}\right)\alpha- \frac{\delta^2 +\delta -1}{\delta^3}
<0,
\ea
\right\}
\ee
for any $\delta\in]\frac{\sqrt{5}-1}{2},1[$. Thus, there exists an integer $N>2,$ such that for any $n>N$,
\be\label{ineq_all}
\left.
\ba{r}
\left(\frac{\varepsilon_1^2+\varepsilon_2+1}{\varepsilon_1}\right) \alpha-\frac{\lambda_n}{\delta \lambda_{n-1} }<0,\\
\frac{1}{\varepsilon_1}\left(1+\frac{1}{\varepsilon_2}\right)\alpha -\frac{\lambda_{n}}{\delta \lambda_{n-1} }<0,\\
\frac{1}{\varepsilon_1}\left(1+\frac{1}{\varepsilon_2}\right)\alpha-\frac{\lambda_n}{\delta \lambda_{n-1} }+\frac{1}{\delta^2}\left(\frac{1}{\delta}-1 \right)<0.
\ea
\right\}
\ee
By $x_{n+1}-x_n=\frac{y_{n+1}-x_{n+1}}{\delta}$, Remark \ref{rem_lower_bound} and Lemma \ref{lem4}, for any $\varepsilon_1,\varepsilon_2>0$ and $M>N+1$, we have
\be\label{ineq4}
a_{M+1}-a_{N+1}&=&\sum\limits_{n=N+1}^{M}(a_{n+1}-a_{n})\\
&\leq& \sum\limits_{n=N+1}^{M}\left[\left(\varepsilon_1+\frac{1+\varepsilon_2}{\varepsilon_1}\right) \alpha-\frac{\lambda_n}{\delta \lambda_{n-1} }\right]\|x_{n+1}-y_n\|^2\nonumber\\
&&+\sum\limits_{n=N+2}^{M}\left[\frac{1}{\varepsilon_1}\left(1+\frac{1}{\varepsilon_2}\right)\alpha-\frac{\lambda_n}{\delta \lambda_{n-1} }+\frac{1}{\delta^2}\left(\frac{\lambda_{n-1}}{\delta \lambda_{n-2}}-1 \right)\right]\|x_n-y_n\|^2\nonumber\\
&&+ \left[\frac{1}{\varepsilon_1}\left(1+\frac{1}{\varepsilon_2}\right)\alpha -\frac{\lambda_{N+1}}{\delta \lambda_{N} }\right]\|x_{N+1}-y_{N+1}\|^2+\xi_M\nonumber\\
&\leq&\xi_M.
\ee
where $\xi_M=\frac{1}{\delta^2}\left(\frac{\lambda_{M}}{\delta \lambda_{M-1}}-1 \right)\|x_{M+1}-y_{M+1}\|^2<+\infty $ from Remark \ref{rem_lower_bound} and Correction step. This together with $a_{n}\geq0$ in (\ref{a_n2}) implies that $\{a_{n}\}_{n\in\bN}$ is  bounded and
\ben
\left.
\ba{r}
0\leq-\sum\limits_{n=N+1}^{\infty}\left[\left(\varepsilon_1+\frac{1+\varepsilon_2}{\varepsilon_1}\right) \alpha-\frac{\lambda_n}{\delta \lambda_{n-1} }\right]\|x_{n+1}-y_n\|^2<+\infty\\
0\leq-\sum\limits_{n=N+2}^{\infty}\left[\frac{1}{\varepsilon_1}\left(1+\frac{1}{\varepsilon_2}\right)\alpha-\frac{\lambda_n}{\delta \lambda_{n-1} }+\frac{1}{\delta^2}\left(\frac{\lambda_{n-1}}{\delta \lambda_{n-2}}-1 \right)\right]\|x_n-y_n\|^2<+\infty,
\ea\right\}
\een
so $\lim\limits_{n\rightarrow\infty}\|x_{n+1}-y_n\|=0$ and $\lim\limits_{n\rightarrow\infty}\|x_n-y_n\|=0.$ By the fact
$\|x_{n+1}-x_n\|=\frac{1}{\delta}\|x_{n+1}-y_{n+1}\|$, we have
$\lim\limits_{n\rightarrow\infty}\|x_{n+1}-x_n\|=0$.

Due to $\|x_n-\bar{x}\|^2\leq a_n$, then $\{x_{n}\}_{n\in\bN}$ is bounded. We can complete the proof by Remark \ref{rem_lower_bound} and the similar methods as in the proof of Theorem \ref{thm21}.
\endproof

\begin{rem}
By the above analysis, it seems that convergence of the proposed algorithm could be still ensured  without  the assumption (A3), but
it is not clear how to prove this as far as we known. Actually, the assumption (A3) is not
restrictive, $g$ is continuous on $\dom g$ when $\dom g$ is an open set (this includes all finite-valued functions) or $g=\delta_C$ for any closed convex set $C$. Moreover,   (A3) holds for any separable lower semicontinuous convex function from \cite[Corollary 9.15]{Bauschke2011Convex}.
\end{rem}

\subsection{Ergodic Convergence Rate for $\delta\in]\frac{\sqrt{5}-1}{2},1]$}
\label{sec_Rate}
Since there are many researches about the convergence rate when $\delta\geq1$,   we  just focus on the case when $\delta\in]\frac{\sqrt{5}-1}{2},1]$.  Actually, the optimal rate of convergence is $\cO(1/n)$ for the extragradient method \cite{rate-convergence}. In this subsection, we investigate the ergodic convergence rate of the sequence $\{y_{n}\}_{n\in\bN}$ for the general case (\ref{primal_problem0}).

From \cite{11.} and \cite[Lemma 2.12]{Proximal-extrapolated}, $x^*\in\cS$ if and only if $x^*\in \dom g$ and
\ben
\max_{x\in \dom g}\Phi(x, x^*) := \langle F(x), x^*-x\rangle + g(x^*) -g(x)=0.
\een
The following theorem shows that  the above criteria can be used to find $x^*$ under a desired accuracy.

\begin{thm}
Let $\{x_{n}\}$ and $\{y_{n}\}$ be generated by Algorithm \ref{algo1}. For any $n_1>N$ and a sufficiently large $J\in\bN$ related to $n_1$, we define
\ben
\hat{\lambda}_j = \sum^{j}_{l={n_1}}\lambda_l +\delta \lambda_{n_1}\quad \textrm{and}\quad \hat{x}_j =\frac{1}{\hat{\lambda}_j}\left( \sum^{j}_{l={n_1}+1}\lambda_l y_l+ (1 +\delta )\lambda_{n_1} x_{n_1}\right)
\een
for any $j>J$, then $\hat{x}_j \in\dom g$ and
\ben
\Phi(x, \hat{x}_j)\leq
\frac{\|x_{n_1}-\bar{x}\|^2+ \delta\lambda_{n_1}\Phi(\bar{x}, x_{n_1-1})+ \frac{1+\varepsilon_2}{\varepsilon_1}\alpha\|x_{n_1}-y_{{n_1}-1}\|^2 }{2\hat{\lambda}_j}, ~~ \forall x\in \cH.
\een
\end{thm}
\proof
First of all, we have by (\ref{rate}) that
\ben
&&2\lambda_n(1 +\delta)	\Phi(\bar{x}, x_n)- 2\lambda_{n} \delta	\Phi(\bar{x}, x_{n-1})\\
&\leq&\|x_n-\bar{x}\|^2-\|x_{n+1}-\bar{x}\|^2\\
&&+\left[\frac{1}{\varepsilon_1}\left(1+\frac{1}{\varepsilon_2}\right)\alpha-\frac{\lambda_n}{\delta \lambda_{n-1}} \right]\|x_n-y_n\|^2  -\left(1-\frac{\lambda_n}{\delta \lambda_{n-1}} \right)\frac{1}{\delta^2}\|x_{n+1}-y_{n+1}\|^2 \\
&&+ \frac{1+\varepsilon_2}{\varepsilon_1}\alpha\|x_n-y_{n-1}\|^2 -\left(\frac{\lambda_n}{\delta \lambda_{n-1} }-\varepsilon_1 \alpha\right)\|x_{n+1}-y_n\|^2.
\een
Since $\|x_n-y_n\|\rightarrow0$ as $ n\rightarrow+\infty$, there exists a sufficiently large $J$ such that for any $j>J$, it holds $\|x_j-y_j\|\leq\|x_{n_1}-y_{n_1}\|\neq0$ (If $\|x_{n_1}-y_{n_1}\|=0$, then $\|x_{n_1+1}-y_{n_1+1}\|\neq0$, else $x_{n_1+1}$ is a solution).
So, we let $\|x_{n_1}-y_{n_1}\|\neq0$ with $n_1>N$. Recalling (\ref{ineq_all})  we deduce for any $j>J$ that
\ben
&&2\left(\lambda_{j}(1 +\delta)	\Phi(\bar{x}, x_{j}) +\sum^
{j-1}_{l=n_1}[\lambda_l (1 + \delta)-\lambda_{l+1}\delta]\Phi(\bar{x}, x_l) \right)\\
&\leq&
\|x_{n_1}-\bar{x}\|^2+\left[\frac{1}{\varepsilon_1}\left(1+\frac{1} {\varepsilon_2}\right)\alpha-\frac{\lambda_{n_1}}{\delta \lambda_{{n_1}-1}} \right]\|x_{n_1}-y_{n_1}\|^2 +\left(\frac{\lambda_{j}}{\delta \lambda_{{j}-1}}-1 \right)\frac{1}{\delta^2}\|x_{{j}+1}-y_{{j}+1}\|^2 \\
&&+ \frac{1+\varepsilon_2}{\varepsilon_1}\alpha\|x_{n_1}-y_{{n_1}-1}\|^2 + \delta\lambda_{n_1}\Phi(\bar{x}, x_{n_1-1})\\
&\leq&
\|x_{n_1}-\bar{x}\|^2+\left[\frac{1}{\varepsilon_1}\left(1+\frac{1} {\varepsilon_2}\right)\alpha-\frac{\lambda_{n_1}}{\delta \lambda_{{n_1}-1}}+\left(\frac{1}{\delta }-1\right) \frac{1}{\delta^2} \right]\|x_{n_1}-y_{n_1}\|^2  \\
&&+ \frac{1+\varepsilon_2}{\varepsilon_1}\alpha\|x_{n_1}-y_{{n_1}-1}\|^2 + \delta\lambda_{n_1}\Phi(\bar{x}, x_{n_1-1})\\
&\leq&
\|x_{n_1}-\bar{x}\|^2+ \frac{1+\varepsilon_2}{\varepsilon_1}\alpha\|x_{n_1}-y_{{n_1}-1}\|^2 + \delta\lambda_{n_1}\Phi(\bar{x}, x_{n_1-1}).
\een
Note that the function $\Phi(\bar{x}, \cdot)$ is convex. Now, applying the Jensen's inequality to the left-hand side of the above inequality and taking
\ben
\lambda_{j}(1 +\delta)+\sum^
{{j}-1}_{l={n_1}}[\lambda_l (1 + \delta)-\lambda_{l+1}\delta]=\sum^{j}_{l={n_1}}\lambda_l +\delta \lambda_{n_1}
\een
into account, we have
\ben
2\left(\sum^{j}_{l={n_1}}\lambda_l +\delta \lambda_{n_1}\right)\Phi(\bar{x}, \hat{x}_{j})\leq\|x_{n_1}-\bar{x}\|^2+ \frac{1+\varepsilon_2}{\varepsilon_1}\alpha\|x_{n_1}-y_{{n_1}-1}\|^2 + \delta\lambda_{n_1}\Phi(\bar{x}, x_{n_1-1}),
\een
where
\ben
\hat{\lambda}_j \hat{x}_j= \lambda_{j}(1 +\delta) x_{j} +\sum^
{j-1}_{l=n_1}[\lambda_l (1 + \delta)-\lambda_{l+1}\delta]x_l=\sum^{j}_{l={n_1}+1}\lambda_l y_l+ (1 +\delta )\lambda_{n_1} x_{n_1}.
\een
Evidently, $\hat{x}_j\in \dom g$ which ends  the proof.
\endproof

Notice that $ \{\lambda_n\}$ has a lower bound $\tau>0$ from Remark \ref{rem_lower_bound}. Fixing $n_1>N$, then we get $\hat{\lambda}_j\rightarrow \infty$ as $j\rightarrow\infty$. This implies $\hat{\lambda}_j\geq (j-n_1)\tau$ and Algorithm \ref{algo1} has the ergodic convergence rate $\cO(1/j)$ when $j>J$.

\subsection{Heuristics on Nonmonotonic Step Sizes}
\label{sec_improved}
Generally speaking, the variable step is more beneficial than a fixed step for the proximal gradient methods. In Algorithm \ref{algo1}, the step size $\{\lambda_n\}_{n\in\bN}$ is updated but in a nonincreasing way, which might be adverse if the algorithm starts in the region with a big curvature of $F$. Namely, the step size in Algorithm \ref{algo1} is  overdependent   on the initial point. For the purpose of obtaining nonmonotonic step sizes, we present an improved algorithm as follows:

\vskip5mm
\hrule\vskip2mm
\begin{algo}
[Improved PEG with nonmonotonic step size.]\label{algo3}
{~}\vskip 1pt {\rm
\begin{description}
\item[{\em Step 0.}] Take $\delta\in]\frac{\sqrt{5}-1}{2},+\infty[ $,  choose  $x_0\in \mathcal{H},$  $\lambda_0>0$, $\gamma\in(0,1)$, $\alpha \in ]0, \kappa(\delta)[$ and a bounded sequence $\{\zeta_n\}$. Set $y_0=x_0$, $x_{1}=\mbox{prox}_{\lambda_0 g}(x_0-\lambda_0 F(x_0))$ and $n=1$.  Choose $\widehat{\lambda}>0$ and a sequence $\{\phi_n\}$ with $\phi_n\in[1,\frac{1+\delta} {\delta}]$ and $\phi_n=1$ when $n\geq n_0$ for given $n_0$.
\item[{\em Step 1.}] \textbf{Prediction:}\\
1.a. Compute
\be
y_{n}&=&x_{n}+\delta (x_{n}-x_{n-1}),\\
\lambda_{n}&=&
         \min~\left\{\phi_{n-1}\lambda_{n-1}, ~~\frac{\alpha\|y_{n}-y_{n-1}\|} {\|F(y_{n})-F(y_{n-1})\|},~~\widehat{\lambda} \right\}.\label{lambda_tilde}
\ee
1.b. Compute \ \ \ \
\ben
 x_{n+1}=\mbox{prox}_{\lambda_n g}(x_n-\lambda_n F(y_n)),
\een
if $x_{n+1}=x_n=y_n$, then stop: $x_{n+1} $ is a solution.
\item[{\em Step 2.}] \blue{\textbf{Correction:}\par
Check
\ben
\|x_{n+1}-x_n\|\leq \zeta_n,
\een
if not hold, set $\lambda_n\leftarrow \gamma \lambda_n$ and return to Step 1.b.}
\item[{\em Step 3.}] Set $n\leftarrow n + 1$ and return to Step 1.\\
  \end{description}
}
\end{algo}
\vskip1mm\hrule\vskip5mm

Since the step size is no longer monotonically decreasing, $a_n\geq0$ in (\ref{a_n1}) is not necessarily valid when $\delta\geq1$, so Algorithm \ref{algo3} implements Correction step for any $\delta\in]\frac{\sqrt{5}-1}{2},+\infty[ $. By $\phi_n\in[1,\frac{1+\delta} {\delta}]$ and $\lambda_{n+1}\leq\phi_n\lambda_{n}$, we can deduce $\delta\lambda_{n+1}\leq(1+\delta)\lambda_{n}$. Then Lemmas \ref{lem0}, \ref{lem1} and \ref{lem4} with (\ref{a_n2}) are still valid for sequences $\{x_n\}$ and $\{y_n\}$ generated by Algorithm \ref{algo3}.

The constant $\widehat{\lambda}$ in Algorithm \ref{algo3} is given only to ensure the upper boundedness of $\{\lambda_n\}$. Hence, it makes sense to choose $\widehat{\lambda}$ quite large. In this case, the step sizes generated are allowed to increase but be bounded from Remark \ref{rem_lower_bound}.
Consequently, it follows from $\phi_n=1$ when $n\geq n_0$ for given $n_0$ that the sequence $\{\lambda_{n}\}_{n>n_0}$ generated by Algorithm \ref{algo3} is monotonically decreasing and then convergent,
\ben
\lim_{n\rightarrow\infty}\lambda_{n}>0,~~ \lim_{n\rightarrow\infty}\frac{\lambda_{n}} {\lambda_{n-1}}=1,
\een
and $\frac{1}{\delta^2}\left(\frac{\lambda_{n}}{\delta \lambda_{n-1}}-1 \right)\|x_{n+1}-y_{n+1}\|^2<+\infty$. Under these conditions, it is not difficult to prove the following convergence theorem by using Lemma \ref{lem4} with (\ref{a_n2}), though we do not know how to choose a proper $n_0$.

\begin{thm}\label{thm3}
Let $\{x_n\}_{n\in\bN}$ be a sequence generated by Algorithm \ref{algo3} with $\delta\in]\frac{\sqrt{5}-1}{2},+\infty[$. Then,  $\{x_n\}_{n\in\bN}$ converges weakly to a solution of problem (\ref{primal_problem0}).
\end{thm}

\section{Further Discussion}
\label{sec_fefurther}
\begin{figure}[htpb]
\centering
\includegraphics[width=0.65\textwidth]{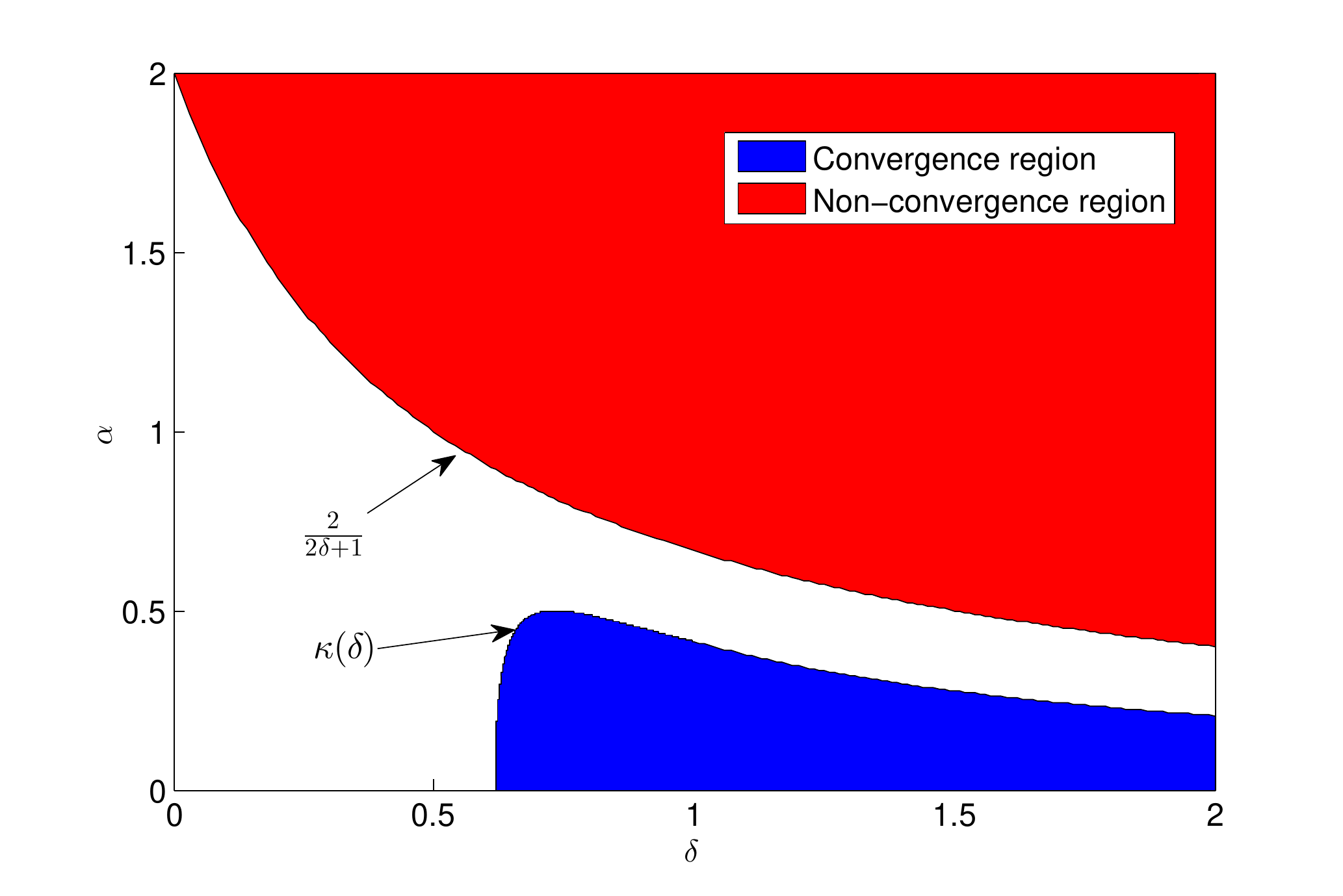}
\caption{The convergence and non-convergence region of the parameters $\delta$ and $\alpha$.}
\label{Fig 2} 
\end{figure}

From the statement above, the condition $\alpha\in]0,\kappa(\delta)[$ for any $\delta\in]\frac{\sqrt{5}-1}{2},+\infty[$ is sufficient to ensure   convergence of the proposed method. In this section, we explain by an extremely simple example that Algorithm \ref{algo1} is not convergent when $\alpha\in]\frac{2}{2\delta+1},+\infty[$ for any $\delta\in]0,+\infty[$. That is to say,  we would derive an upper bound of $\alpha$  to guarantee   the convergence of Algorithm \ref{algo1}, but Algorithm \ref{algo1} with $(\delta,\alpha)$   in some regions remains to be further studied, see Fig. \ref{Fig 2}.

Consider the simplest optimization problem
$$
\min\limits_{x\in\bR^m}\frac{1}{2}\|x\|^2.
$$
Obviously, it can be formulated as a special case of problem (\ref{primal_problem})  with $F=I$ (the identity operator), $L=1$ and $C=\bR^m$. Followed by the updates of Algorithm \ref{algo1}, we have
\be\label{ite_example}
x_{n+1}=(1-\lambda_n-\delta\lambda_n)x_n+\delta\lambda_n x_{n-1},~\delta\in]0,+\infty[.
\ee
For any $\widetilde{d}_n$, $\widehat{d}_n \in \bR$, if
\be\label{weida}
\widetilde{d}_n+\widehat{d}_n=1-\lambda_n-\delta\lambda_n ~~\mbox{and}~~\widetilde{d}_n \widehat{d}_n=-\delta\lambda_n,
\ee
then we can rewrite (\ref{ite_example}) as
\be\label{ite_rewrite}
x_{n+1}-\widetilde{d}_n x_n=\widehat{d}_n(x_n-\widetilde{d}_nx_{n-1}) ~~\mbox{or}~~x_{n+1}-\widehat{d}_nx_n=\widetilde{d}_n(x_n-\widehat{d}_nx_{n-1}).
\ee

By (\ref{weida}) and Vieta's Theorem, we have
$$
\widetilde{d}_n,\widehat{d}_n=\frac{1-\lambda_n-\delta \lambda_n\pm\sqrt{(1-\lambda_n-\delta\lambda_n)^2+4\delta\lambda_n}}{2}.
$$
If $\max\left\{|\widetilde{d}_n|,~|\widehat{d}_n|\right\}>1$, then the iterative (\ref{ite_rewrite}) is not convergent. As a result, (\ref{ite_example}) is not convergent either. Namely, if
\ben
\lambda_n>\frac{2}{2\delta+1},\quad \forall \delta>0,
\een
then the iterative (\ref{ite_example}) is not convergent. By Remark \ref{rem_lower_bound} and $L=1$, the convergence of Algorithm \ref{algo1} can not be guaranteed if $\lambda_0>\frac{2}{2\delta+1}$  and $\alpha\in]\frac{2}{2\delta+1},+\infty[$ for any $\delta\in]0,+\infty[$.

\section{Numerical Experiments}
\label{sec_experiments}
In this section,
we perform Algorithm \ref{algo3} \footnote{All codes are available at http://www.escience.cn/people/changxiaokai/Codes.html} (denoted by ``IPEG")  for solving  some randomly   generated minimization problems over difficult nonlinear constraints. The following state-of-the-art algorithms are compared to    investigate the computational efficiency of IPEG:
\begin{itemize}
\item  Tseng's forward-backward-forward splitting method used as in  \cite[Section 4]{Proximal-extrapolated} (denoted by ``TFBF"), with $\beta= 0.7, \theta= 0.99$;
\item Proximal extrapolated gradient methods \cite[Algorithm 2]{Proximal-extrapolated} (denoted by ``PEG"), with line search and $\alpha= 0.41, \sigma= 0.7$;
\item Modified projected gradient method \cite{yang} (denoted by ``MPG"), with $\alpha= 0.41, \delta= 1.01$.
\item FISTA \cite{Nesterov1983} with standard linesearch (denoted by ``FISTA"), with $\beta= 0.7, \lambda_0= 1$;
\end{itemize}

We denote the random number generator by $seed$ for generating data again in Python 3.8. All experiments are performed on an Intel(R) Core(TM) i5-4590 CPU@ 3.30 GHz PC with 8GB of RAM running on 64-bit Windows operating system.\par

Since solutions of (\ref{primal_problem0}) coincide with zeros of the residual function
\ben
r(x, y):=\|y-\mbox{prox}_{\lambda_n g} (x-\lambda F(y))\|+\|x-y\|,
\een
for some positive number $\lambda$, and $r_n:=r(x_n, y_n)=\|x_{n+1}-y_n\|+\|x_n-y_n\| =0$ implies $x_{n+1}=x_n=y_n$, thus we use $r_n<\epsilon$ with given $\epsilon=10^{-6}$ to terminate our algorithms, and the same $\epsilon$ is used to terminate PEG, MPG, FB and FISTA. In particular for TFBF, we use
\ben
r_n:=\|x_n-\mbox{prox}_{\lambda_n g} (x_n-\lambda F(x_n))\|\leq\epsilon
\een
as in \cite{Proximal-extrapolated}.

We generate $\lambda_0$ as in \cite{M-extra},  choose $y_{-1}$ as a small perturbation
of $y_0$ and take $\lambda_0 = \frac{\|y_{-1}-y_0\|}{\|F(y_{-1})-F(y_0)\|}$. This gives us an approximation of the local inverse Lipschitz constant of $F$ at $y_0$. There are many choices of the sequence $\{\phi_n\}_{n\in\bN}$, but in the earlier iterations the large range of $\lambda_n$ is benefit for selecting proper step size, we thus use
\be\label{phi0}
\phi_n=\left\{
\ba{rl}\frac{1+\delta}{\delta},~&~\mbox{if}~n\leq\hat{n};\\
\frac{1+\delta+n-\widehat{n}}{\delta+n-\widehat{n}},~&~\mbox{if}~n>\widehat{n},
\ea
\right.
\ee
for a given $\hat{n}\in\bN$. In this section, we fix  $\hat{n}=500$ and $n_0=1000$. For applying Correction step, we use $\gamma=0.7$ and (\ref{zeta}) with $\zeta_{\min}=10^{-6}$ and $\mu=\nu=10$.

We report the number of iterations (Iter), the number of
proximal operators ($\#$ prox), the number of $F$ ($\# F$) and the computing time (Time) measured in seconds. Note that the number of iterations equals that of  proximal operators for PEG and IPEG, and is 2 smaller than that of $F$ for IPEG, we thus report the number of iterations and the number of $F$ for PEG and only the number of iterations for IPEG. The bold letter indicates the best results in the following tables.

\begin{prob}\label{pro_2}
The first  problem (called Sun's problem) was considered in \cite{9.,15.,yang}, and the Lipschitz-continuous and monotone operator was generated by
\ben
F(x)=G(x)+H(x),
\een
where
\ben
G(x)&=&(g_1(x),g_2(x),\ldots,g_m(x)),\\
g_i(x)&=&x_{i-1}^2+x_{i}^2+x_{i-1}x_i+x_{i}x_{i+1}, \ \ i=1,2,\ldots,m, \ \ x_0=x_{m+1}=0.
\een
and $ H(x)=Ex+c.$ Here $E$ is a square matrix $m\times m$ defined by
\ben
e_{ij}= \begin{cases}
         4 , & \ \ j=i  \\
          1,& \ \ i-j=1  \\
          -2,& \ \ i-j=-1  \\
          0, & otherwise, \end{cases}
\een
and $c=(-1,-1,\ldots,-1).$ We choose the feasible set $C$ as $C_1=\bR_+^m$ and $C_2 = \{x\in\bR^m_+~|~ \sum_{i=1}^m x_i =m\}$.
\end{prob}

\begin{figure}[htpb]
\centering
\subfigure[$r_{n}$.]{
\includegraphics[width=0.45\textwidth]{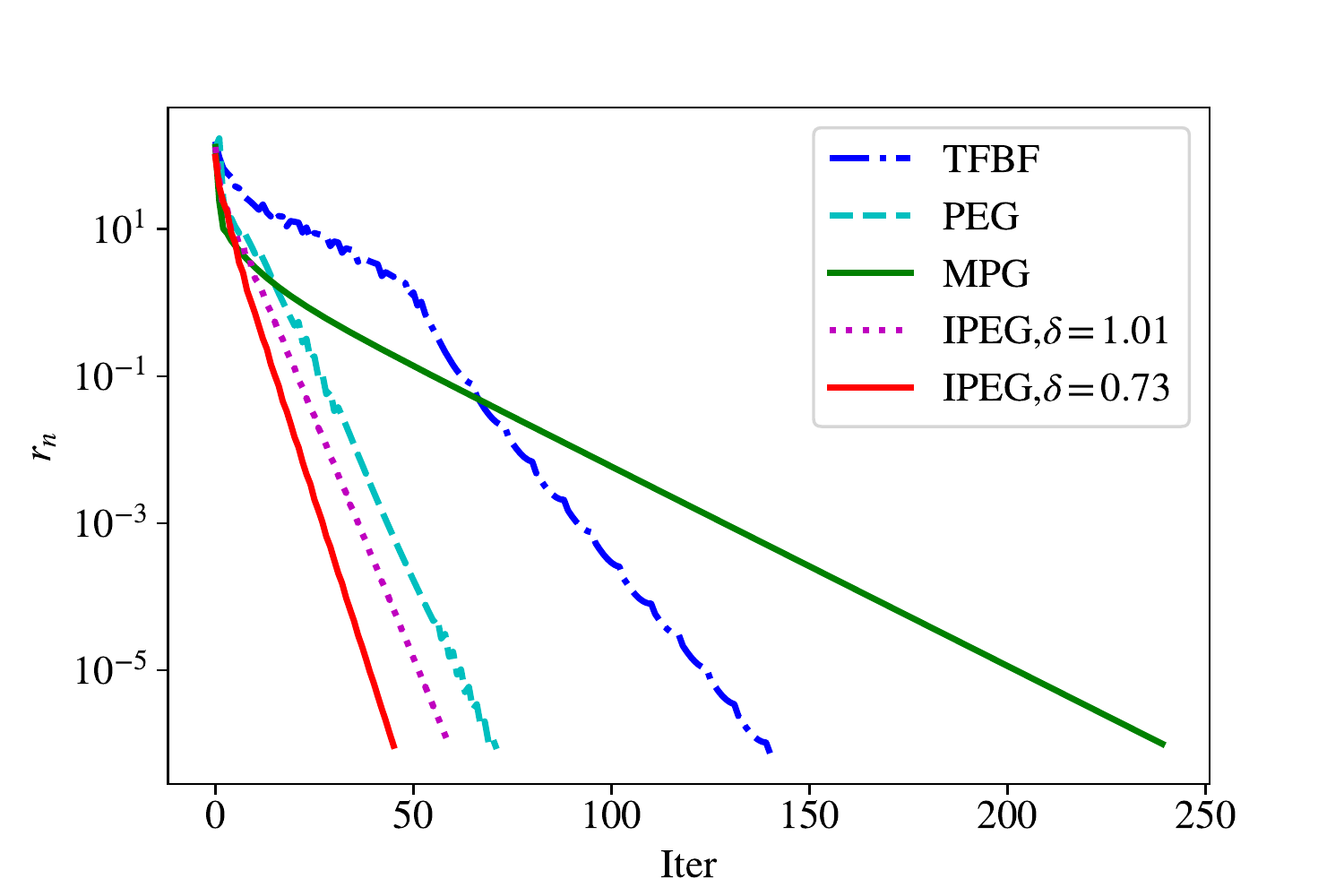}}
\subfigure[$\lambda_n$.]{
\includegraphics[width=0.45\textwidth]{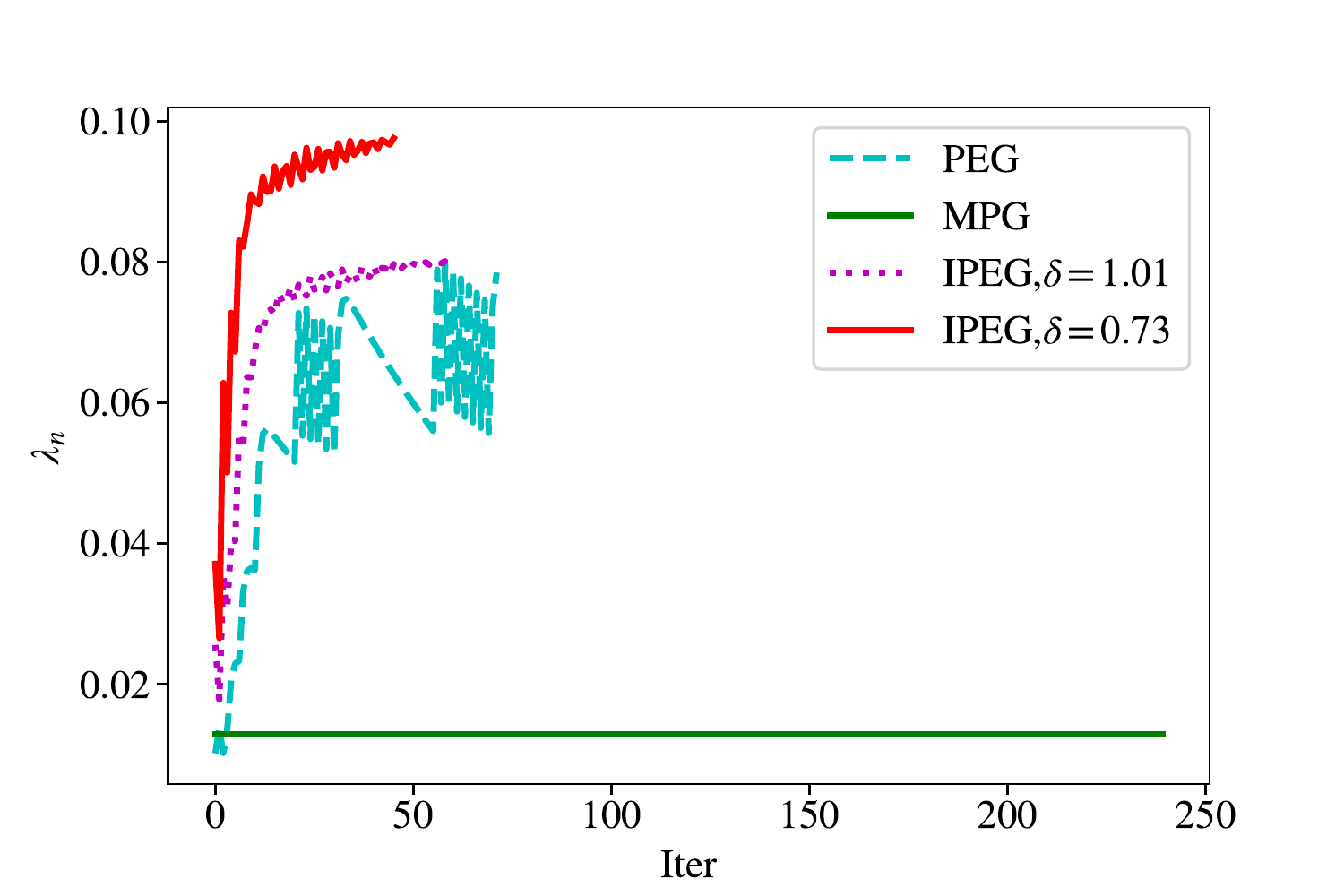}}
\caption{Comparison of $r_{n}$ and $\lambda_n$ for solving Problem \ref{pro_2} with $C=C_1$, $d=10^3$. }
\label{Fig 3} 
\end{figure}

For Problem \ref{pro_2}, the initial point $x_0$ is generated uniformly randomly from $[-10, 10]^d$. For every $d = 10^3, 10^4, 10^5$ and every $C$ above, the test results are listed in Table \ref{table2}. In addition, we show the evolutions of $r_{n}$ and $\lambda_n$ with respect to Iter for solving Problem \ref{pro_2} with $C=C_1$, $d=10^3$ in Fig. \ref{Fig 3}.

\begin{table}[ht]
\caption{ Results for Problem \ref{pro_2} with different $d$ and $C$.}\label{table2}
\vspace{.1 in}
\center
\small
\begin{tabular}{|c|c|cccc|ccc|cc|cc|cc|}
\hline
\multirow{3}{0.5cm}{$C$}&\multirow{3}{0.5cm}{$d$}& \multicolumn{4}{|c|}{TFBF}&\multicolumn{3}{|c|}{PEG} & \multicolumn{2}{|c|}{MPG}&\multicolumn{2}{|c|}{IPEG} &\multicolumn{2}{|c|}{IPEG}  \\
&&&&&&&&&&&\multicolumn{2}{|c|}{($\delta=1.01$)} &\multicolumn{2}{|c|}{($\delta=0.73$)}\\
 & & Iter & $\#$ prox& $\#$ F & Time& Iter &$\#$ F & Time& Iter & Time& Iter & Time& Iter & Time\\ \hline
\multirow{3}{0.5cm}{$C_1$}&	$10^3$	&	141	&	294	&	435	&	0.05	&	73	&	143	&	0.02	&	 243	&	0.03	&	62	&	0.01	&	\textbf{48}	&	0.01	\\
&	$10^4$	&	163	&	341	&	504	&	0.1	&	76	&	149	&	0.04	&	262	&	0.03	&	66	&	 0.02	&	\textbf{50}	&	0.01	\\
&	$10^5$	&	174	&	365	&	539	&	2.21	&	80	&	157	&	0.76	&	284	&	1.23	&	70	&	 0.32	&	\textbf{53}	&	0.31	\\ \hline
\multirow{3}{0.5cm}{$C_2$}&	$10^3$	&	139	&	292	&	431	&	0.05	&	78	&	154	&	0.02	&	 229	&	0.03	&	77	&	0.01	&	\textbf{63}	&	0.01	\\
&	$10^4$	&	145	&	305	&	450	&	0.31	&	83	&	164	&	0.09	&	249	&	0.24	&	83	&	 0.08	&	\textbf{67}	&	0.07	\\
&	$10^5$	&	170	&	359	&	529	&	4.89	&	88	&	174	&	1.44	&	270	&	3.41	&	88	&	 1.13	&	\textbf{71}	&	1.01	\\ \hline
\end{tabular}
\end{table}

\begin{prob}\label{pro_3}
The second test problem is the so-called Kojima-Shindo Nonlinear Complementarity Problem (NCP),  considered in \cite{10.,P-G}, where $m=4$ and the mapping $F$ is defined by
\ben
F(x_1, x_2, x_3, x_4) =\left(\ba{l}
3x_1^2+ 2x_1x_2 + 2x_2^2 + x_3 + 3x_4-6\\
2x_1^2 + x_1 + x_2^2 + 10x_3 + 2x_4-2\\
3x_1^2 + x_1x_2 + 2x_2^2 + 2x_3 + 9x_4-9\\
x_1^2 + 3x_2^2 + 2x_3 + 3x_4-3
\ea
\right).
\een
The feasible set is $C = \{x\in\bR^4_+~|~ x_1+x_2+x_3+x_4 =4\}$ and $g(x)=l_{C}(x)$.
\end{prob}

\begin{figure}[htbp]
\centering
\subfigure[$r_{n}$.]{
\includegraphics[width=0.45\textwidth]{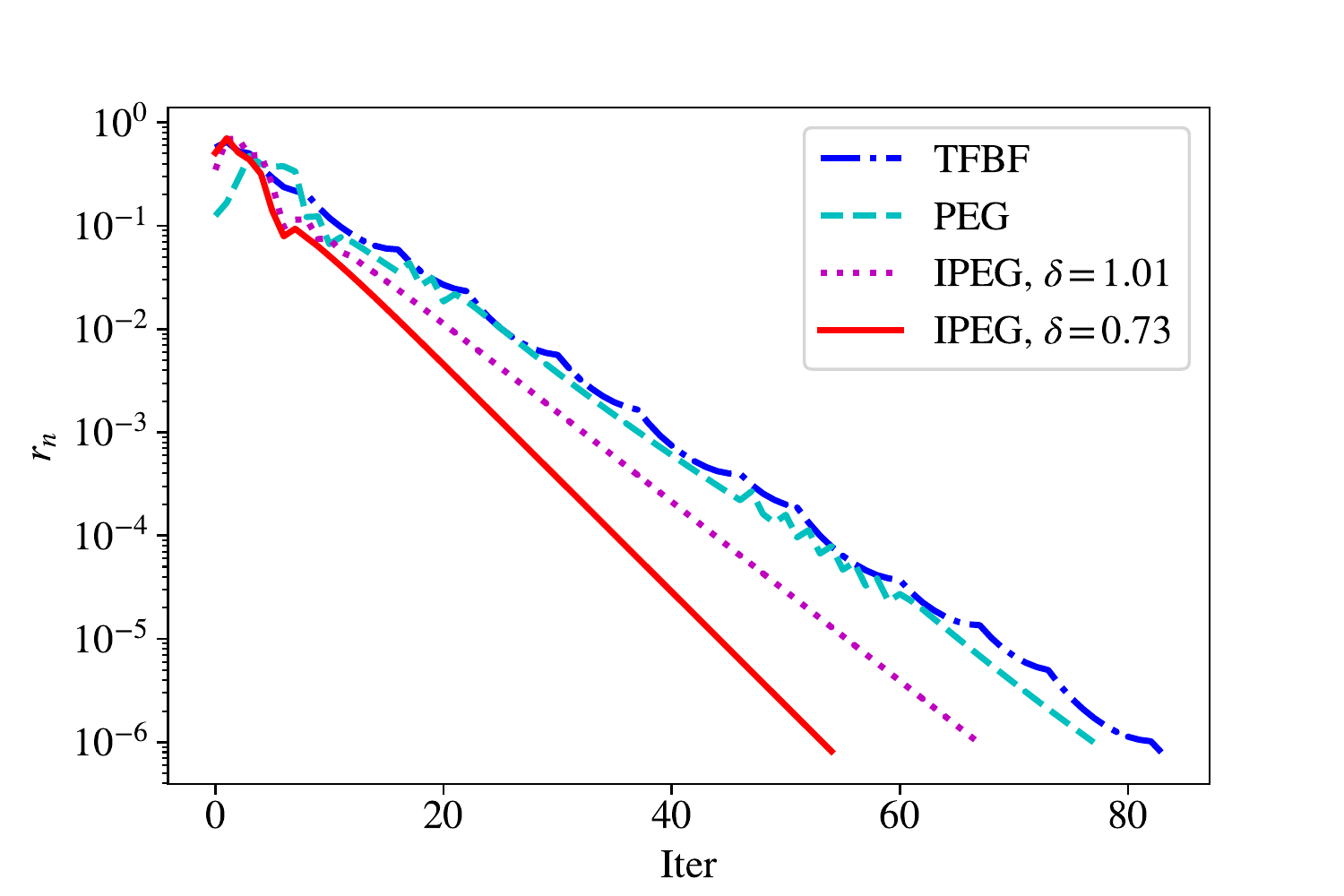}}
\subfigure[$\lambda_n$.]{
\includegraphics[width=0.45\textwidth]{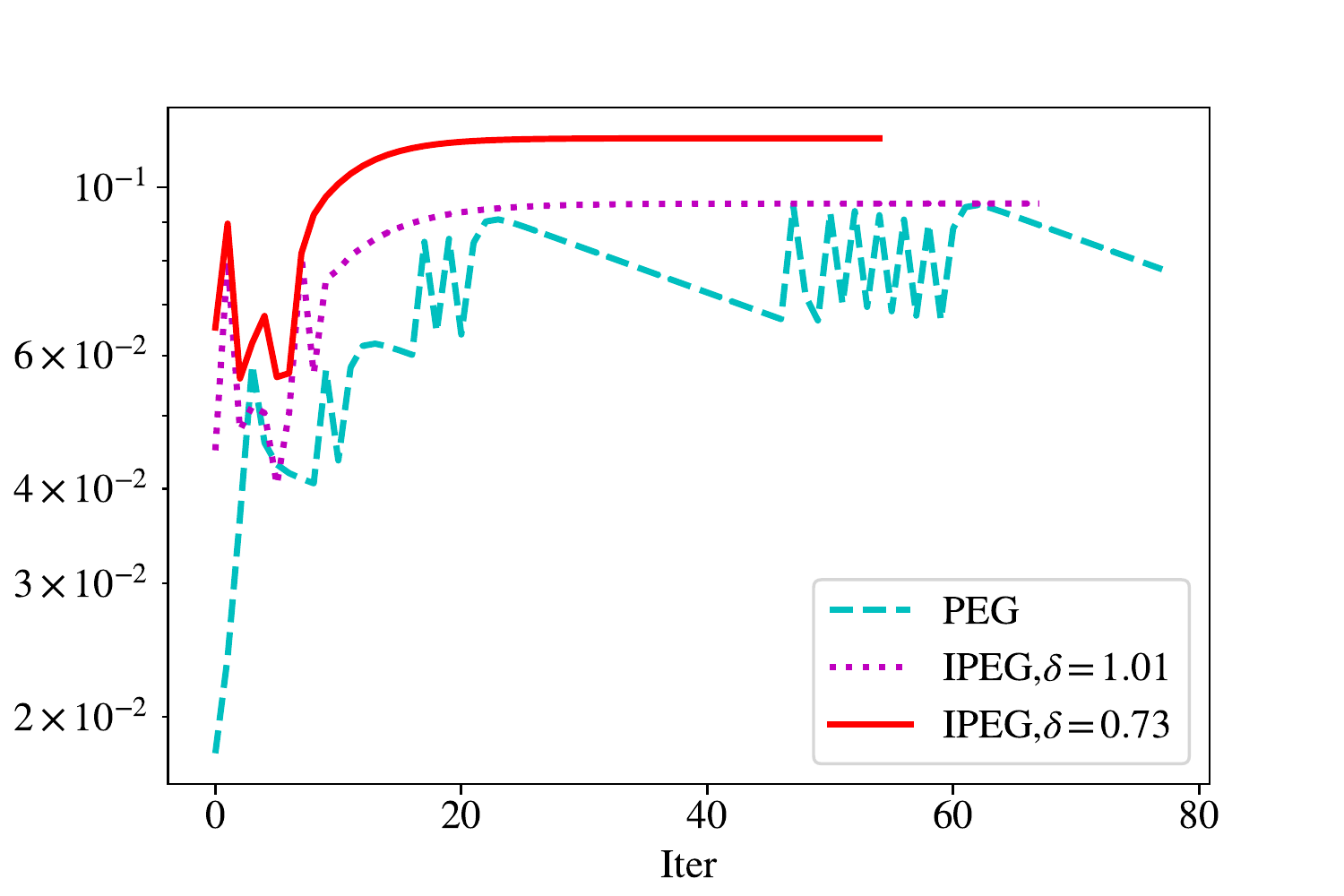}}
\caption{Comparison of $r_{n}$ and $\lambda_n$ for solving Problem \ref{pro_3} with $x_0=(1,1,1,1)$. }
\label{Fig 4} 
\end{figure}

We choose three particular starting points: $(0, 0, 0, 0)$, $(1, 1, 1, 1)$ and $(0.5, 0.5, 2, 1)$. The numerical results are reported in Table \ref{table3} and the evolutions of $r_{n}$ and $\lambda_n$ with respect to Iter for solving Problem \ref{pro_2} with $x_0=(1,1,1,1)$ are shown in Fig. \ref{Fig 4}.

\begin{table}[ht]
\caption{ Results for Problem \ref{pro_3} with different $x_0$.}\label{table3}
\vspace{.1 in}
\center
\small
\begin{tabular}{|c|cccc|ccc|cc|cc|}
\hline
\multirow{2}{0.5cm}{$x_0$}& \multicolumn{4}{|c|}{TFBF}&\multicolumn{3}{|c|}{PEG} &\multicolumn{2}{|c|}{IPEG($\delta=1.01$)}&\multicolumn{2}{|c|}{IPEG($\delta=0.73$)}  \\
& Iter & $\#$ prox& $\#$ F & Time& Iter &$\#$ F & Time& Iter & Time& Iter & Time\\ \hline
 $(0,0,0,0)$	&	81	&	173	&	254	&	0.02	&	82	&	164	&	0.1	&	72	&	0.01	&	 \textbf{58}	&	0.01	\\
$(1,1,1,1)$	&	84	&	177	&	261	&	0.02	&	79	&	156	&	0.1	&	70	&	0.01	&	 \textbf{56}	&	0.01	\\
$(0.5,0.5,2,1)$	&	88	&	186	&	274	&	0.02	&	85	&	169	&	0.1	&	75	&	0.01	&	 \textbf{59}	&	0.01	\\ \hline
\end{tabular}
\end{table}

\begin{prob}\label{pro_4}
The third problem is HpHard problem, considered as in \cite{yang,Proximal-extrapolated}. Let $F(x) = Mx +q$ with $M = NN^T + S + D$ and $q \in \bR^m$, where $N$, $D$ and $S\in\bR^{m\times m}$, $S$ is a skew-symmetric matrix, every entry of $N$ and $S$ is uniformly generated from $(-5, 5)$. The matrix $D$ is diagonal and its diagonal entry is uniformly generated from $(0, 0.3)$. Every entry of $q$ is uniformly generated from $(-500, 0)$. The feasible set is $C = \{x\in\bR^m_+~|~ \sum_{i=1}^mx_i=m\}$ and $g(x)=l_{C}(x)$.
\end{prob}

\begin{figure}[htp]
\centering
\subfigure[$r_{n}$.]{
\includegraphics[width=0.45\textwidth]{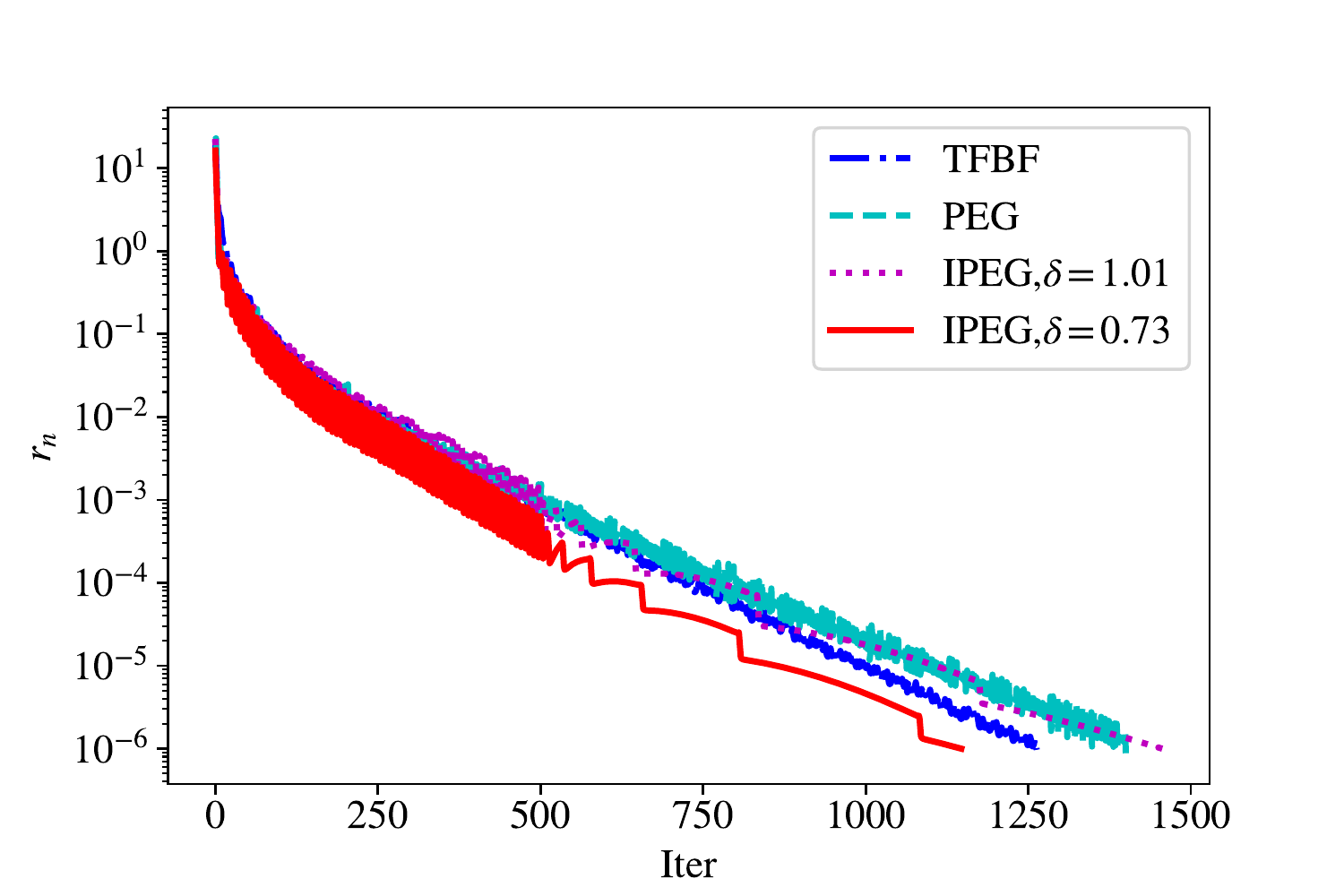}}
\subfigure[$\lambda_n$.]{
\includegraphics[width=0.45\textwidth]{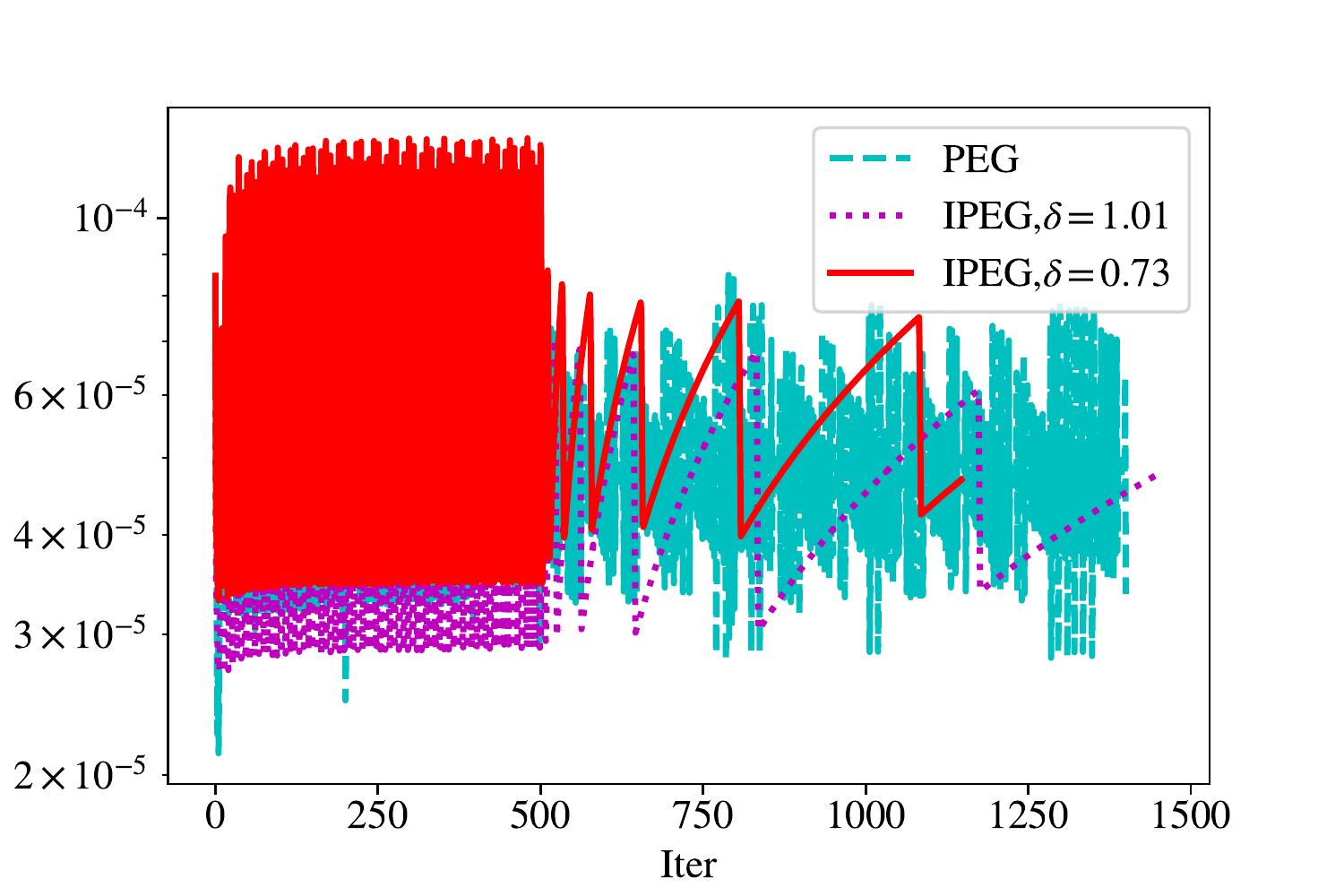}}
\caption{Comparison of $r_{n}$ and $\lambda_n$ for solving Problem \ref{pro_4} with $seed=1$ and $d=500$.  }
\label{Fig 5} 
\end{figure}

\begin{table}[htpb]
\caption{ Results for Problem \ref{pro_4} with different cases.}\label{table4}
\vspace{.1 in}
\center
\small
\begin{tabular}{|c|c|cccc|cc|cc|cc|}
\hline
\multirow{2}{0.5cm}{$seed$}&\multirow{2}{0.5cm}{$m$}& \multicolumn{4}{|c|}{TFBF}&\multicolumn{2}{|c|}{PEG} & \multicolumn{2}{|c|}{IPEG($\delta=1.01$)}&\multicolumn{2}{|c|}{IPEG($\delta=0.73$)}  \\
 & & Iter & $\#$ prox& $\#$ F & Time& Iter  & Time& Iter & Time& Iter & Time\\ \hline
\multirow{3}{0.5cm}{1}	&	500	&	1066	&	2279	&	3345	&	0.25	&	1185	&	0.17	&	 1185	&	0.13	&	\textbf{972}	&	0.09	\\
&	1000	&	1155	&	2469	&	3624	&	1.75	&	1323	&	0.93	&	1268	&	0.42	 &	\textbf{1033}	&	0.39	\\
&	5000	&	1389	&	2969	&	4358	&	56.87	&	1575	&	27.89	&	1630	&	26.43	 &	\textbf{1326}	&	22.86	\\
 \hline\hline
\multirow{3}{0.5cm}{2}	&	500	&	1270	&	2715	&	3985	&	0.29	&	1447	&	0.19	&	 1480	&	0.14	&	\textbf{1165}	&	0.12	\\
&	1000	&	1134	&	2424	&	3558	&	1.56	&	1274	&	0.86	&	1262	&	0.41	 &	\textbf{1028}	&	0.40	\\
&	5000	&	1365	&	2918	&	4283	&	55.74	&	1554	&	33.91	&	1603	&	29.94	 &	\textbf{1303}	&	25.64	\\ \hline
\end{tabular}
\end{table}

For every $m$, as shown in Table \ref{table4}, we have generated randomly two different $M$ and $q$ with $seed=1$ and $2$. For all tests, we take $x_0=(1,1,\cdots,1)$. Since $F$ is an affine operator, the number of iterations is 2 smaller than that of $F$ for PEG, thus we  just report  the number of iterations.

\begin{prob}\label{pro_1}
The fourth example is a sparse logistic regression problem for binary classification. Let $(h_i,l_i)\in\bR^n \times \{\pm1\}, i = 1,\cdots,m$ be the training set, where $h_i\in\bR^n$ is the feature vector of each data
sample, and $l_i$ is the binary label. The formulation of sparse logistic regression reads
\be\label{SLR}
\min\limits_{x\in\bR^n}\phi(x):=\mu\|x\|_1 + \frac{1}{m}\sum_{i=1}^m \log
(1 + e^{l_ih^T_ix}),
\ee
where $\mu>0$ and is set to be $0.005\|H^Tl\|_{\infty}$ in the numerical test.
\end{prob}

Let $K_{ij}=-l_ih_{ij}$ and set $\hat{f}(y)=\sum^m_{i=1} \log(1+\exp(y_i))$. Then the objective in (\ref{SLR}) is $\phi(x)= f(x)+g(x)$ with $g(x) =\mu \|x\|_1$ and $f(x)=\hat{f}(Kx)$. It is easy to derive that $L_{\nabla \hat{f}}=\frac{1}{4}$. Thus, $L_{\nabla f}=\frac{1}{4}\|K^TK\|$. We take three popular datasets from LIBSVM \footnote{https://www.csie.ntu.edu.tw/$\sim$cjlin/libsvmtools/datasets/}: \textbf{w7a} with $m=24692$, $n=300$, \textbf{a9a} with $m = 32561$, $n = 123$ and \textbf{real-sim} with $m = 72309$, $n = 20958$.

\begin{figure}[htp]
\centering
\subfigure[$r_{n}$.]{
\includegraphics[width=0.45\textwidth]{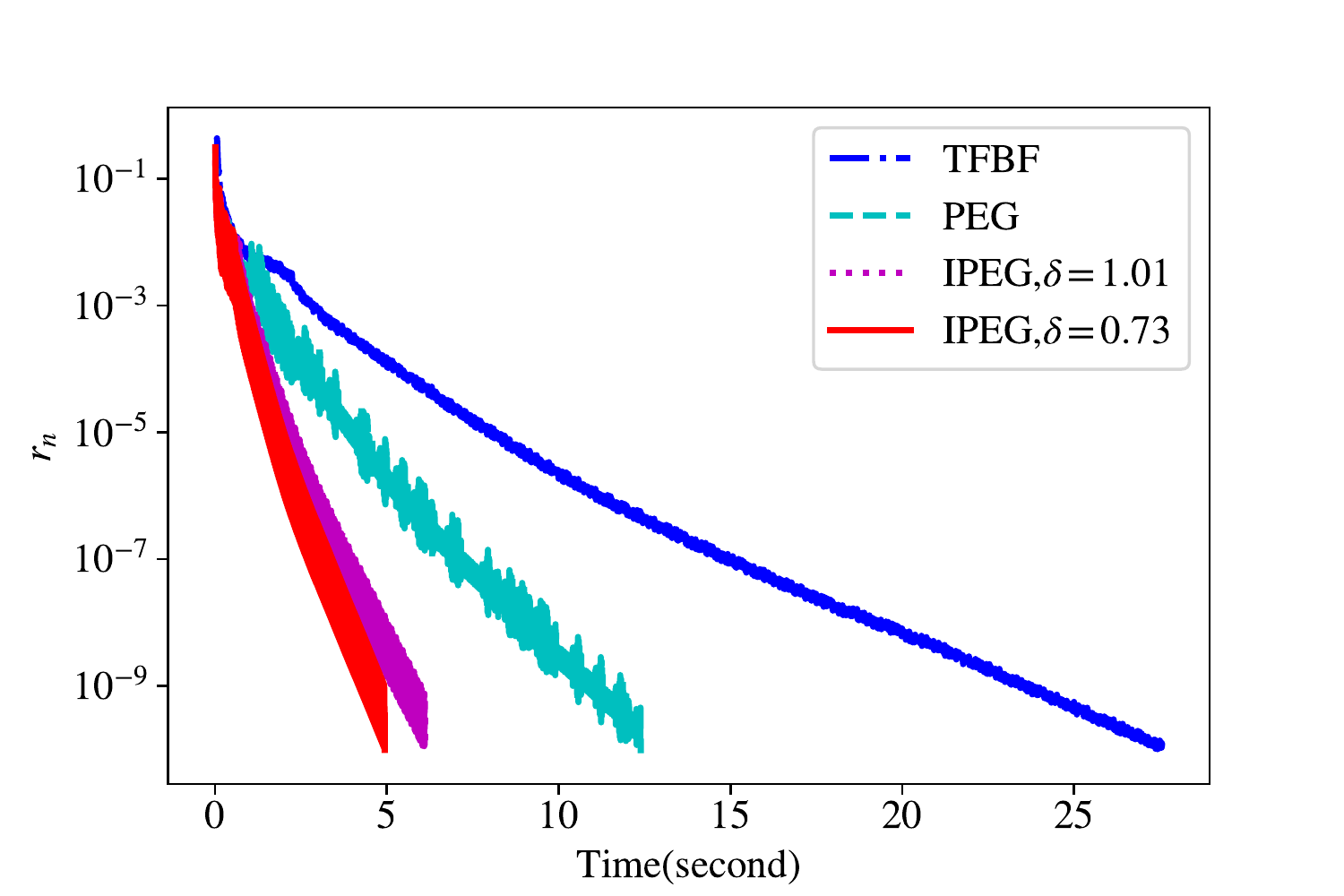}}
\subfigure[$\lambda_n$.]{
\includegraphics[width=0.45\textwidth]{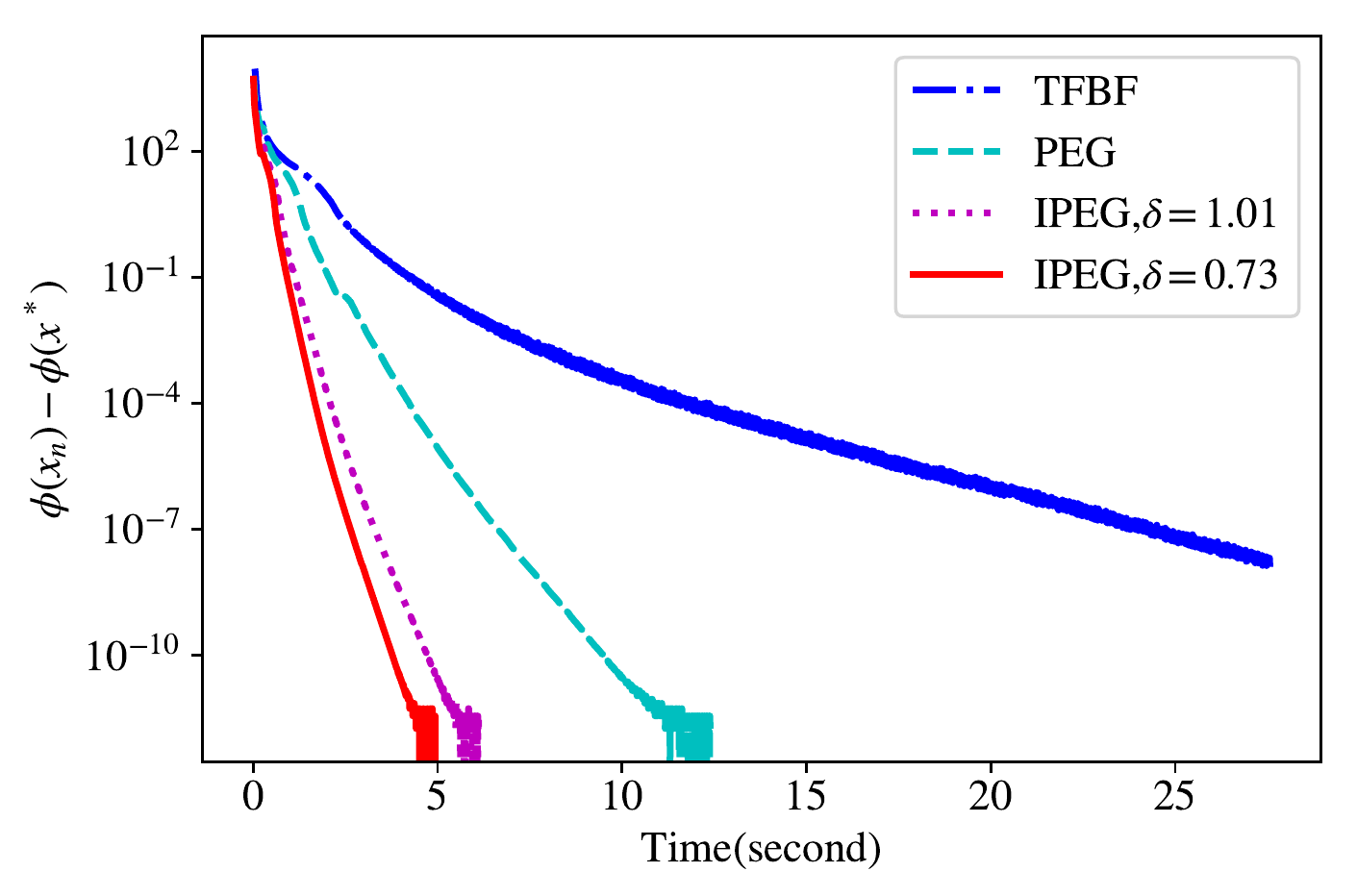}}
\caption{Comparison of $r_{n}$ and $\phi(x_n)-\phi(x^*)$ for solving Problem \ref{pro_1} with data ``a9a".  }
\label{Fig 6} 
\end{figure}

Since $f$ is convex and $F=\nabla f$, we apply IPEG to (\ref{SLR}) without Correction step. We use $\epsilon=10^{-10}$ to terminate all the algorithms for getting more accurate solution, and choose the smallest objective value among all methods and set it to $\phi(x^*)$. The results are shown in Table \ref{table5}. To illustrate how does the value  $\phi(x_n)-\phi(x^*)$ and $r_n$ change over times, we give two convergence plots for data ``a9a" in Fig. \ref{Fig 6}.

\begin{table}[htpb]
\caption{ Results for Problem \ref{pro_1}.}\label{table5}
\vspace{.1 in}
\center
\small
\begin{tabular}{|c|cccc|ccc|cc|cc|}
\hline
\multirow{2}{0.5cm}{data}& \multicolumn{4}{|c|}{TFBF}&\multicolumn{3}{|c|}{PEG} & \multicolumn{2}{|c|}{IPEG($\delta=1.01$)}&\multicolumn{2}{|c|}{IPEG($\delta=0.73$)}  \\
 &  Iter & $\#$ prox& $\#$ F & Time& Iter & $\#$ F & Time& Iter & Time& Iter & Time\\ \hline
w7a &		971	&	1950	&	2867	&	4.1	&	968	&1933&	2.9	&	827	&	1.6	 &	\textbf{716}	&	1.4	\\
a9a	&		6758	&	14439	&	21197	&	27.8	&	4241	&8601&	12.2	&	 3498	&	6.1	&	\textbf{2844}	&	5.0	\\
real-sim&		3984	&	8510	&	12494	&	153.8	&	2651	&	5312	&	70.9&2230	&	35.1	 &	\textbf{1796}	&	32.8	\\
\hline
\end{tabular}
\end{table}

To summarize our numerical experiments on Problems \ref{pro_2}-\ref{pro_1}, we want to make some observations. Firstly, the advantage of IPEG in comparison with other algorithms is a larger interval for possible step size $\lambda_n$, see Fig. \ref{Fig 3}(b), Fig. \ref{Fig 4}(b) and Fig. \ref{Fig 5}(b), which resulted from the proper choice of $\delta$ and the larger value of $\alpha$.

Secondly, we observed that for the majority of the test problems, IPEG is more efficient than other algorithms in both the number of iterations and the CPU time. Furthermore, IPEG with $\delta= 0.73$ performs efficiently than that with $\delta= 1.01$ from the convergence plots of $r_n$ shown in Fig. \ref{Fig 3}(a), Fig. \ref{Fig 4}(a) and Fig. \ref{Fig 5}(a),  which is extremely due to the larger step size $\lambda_n$ and the use of only one value of the mapping required per iteration. Although linesearch is involved in Correction step, the condition required is so weak that the linesearch is not started for many problems.

In addition, since MPG \cite{yang} adopted nonincreasing step sizes, it is adverse when starting in the region with a big curvature of $F$, see Fig. \ref{Fig 3}(b) and the results of MPG for Problem \ref{pro_2}. From Fig. \ref{Fig 5}, the step sizes generated by IEPG have fluctuated within a range at the first 500 iterations, after that the range decreases as we use (\ref{phi0}) with $\widehat{n}=500$ to control the increase of step sizes.

\section{Conclusions}
\label{sec_conclusion}
Without the knowledge of Lipschitz constant, we have proposed a proximal extrapolated gradient method using a prediction-correction procedure to determine stepsizes, and improved it numerically with non-monotonic step size.  The method extended the range of parameters (considering the case of $\delta<1$ ) and obtained a larger step size than the existing methods by using correction step. Finally, a number of experiments illustrate that the proposed method is efficient, and the improvement can be resulted from the larger step size.

In addition, we have shown by an extremely simple example that our method is not convergent if $\lambda_0,~\alpha\in]\frac{2}{2\delta+1},+\infty[$ for any $\delta>0$. From Fig. \ref{Fig 3}, the convergence of the proposed method remains unknown for $(\delta,\alpha)$ in some regions. Especially for $\delta\in]0,\frac{\sqrt{5}-1}{2}]$, it remains to be explored whether there are any (larger) $\alpha>0$ such that Algorithms \ref{algo1} and \ref{algo3} are convergent. Perhaps our method without the correction step is convergent as well, and can be generalized to other methods that need to estimate the Lipschitz constant. We leave this as an interesting topic for our future research.

\begin{acknowledgements}
The research of Xiaokai Chang was supported by  the Hongliu Foundation of First-class Disciplines of Lanzhou University of Technology. The project was supported by  the National Natural Science Foundation of China under Grant 61877046 and the Natural Science Basic Research Plan in Shaanxi Province of China (2017JM1014).
\end{acknowledgements}

\end{document}